\documentclass[12pt]{amsart}
\usepackage{amssymb,amscd}
\usepackage[dvipdfmx]{graphicx,color}

\usepackage{graphicx}

\usepackage{amsmath}

\newtheorem*{Whitney towers}{Theorem~\ref{Whitney towers}}
\newtheorem*{h-towers}{Theorems ~\ref{half} \& \ref{$(n)$-solvable}}

\newtheorem*{surgery curves}{Theorem~\ref{surgery curves}}
\newtheorem*{cg=0}{Theorem~\ref{vanish}}

\newtheorem{thm}{Theorem}[section]
\newtheorem{mth}[thm]{Main Theorem}
\newtheorem{pr}[thm]{Proposition} 
\newtheorem{fact}[thm]{Fact} 
\newtheorem{lem}[thm]{Lemma}
\newtheorem{cor}[thm]{Corollary}

\newtheorem{cla}[thm]{Claim}

\theoremstyle{definition}
\newtheorem{defn}[thm]{Definition}
\newtheorem{que}[thm]{Question}
\newtheorem{prob}[thm]{Problem}
\newtheorem{note}[thm]{Note}

\numberwithin{equation}{section}
\numberwithin{figure}{section}

\newcommand{\x}{\times}
\newcommand{\np}{\newpage}
\newcommand{\Z}{\mathbb{Z}}
\newcommand{\N}{\mathbb{N}}

\newcommand{\R}{\mathbb{R}}

\def\yen{{\setbox0=\hbox{Y}Y\kern-.97\wd0\vbox{hrule height.lex width.98%
\wd0\kern.33ex\hrule height.lex width.98\wd0\kern.45ex}}}

\def\np{\newpage}

\begin{document}
\pagestyle{plain}

\title{
A new way to detect pairs of non-cobordant surface-links
which 
Orr invariant, Cochran sequence, Sato-Levine invariant, 
and alinking number 
cannot distinguish
}

\author{Eiji Ogasa}


\begin{abstract}
We submit a new way to detect pairs of non-cobordant surface-links.  
We find a new example of a pair of non-cobordant surface-links with the following properties: 
Orr invariant, Cochran sequence, Sato-Levine invariant, 
the alinking number and 
one of Stallings's theorems  
cannot distinguish them. However our new way can distinguish them. 
\end{abstract} 

\thanks{\hskip-4mm E-mail: 
ogasa@mail1.meijigakuin.ac.jp \newline
}

\date{}

\maketitle

\tableofcontents

\section{Introduction and the main result}\label{Introduction}


\noindent 
We work in the smooth category. 
Let $\mu\in\N$.
Let $i\in\{1,...,\mu\}$. 
Let $F_i$ be a connected closed oriented surface. 
A {\it surface-$(F_1,...,F_\mu)$-link}, or  {\it $(F_1,...,F_\mu)$-link},   
is a submanifold $L=(K_1,...,K_\mu)$ of $S^4$ such that $K_i$ is  diffeomorphic to $F_i$ 
for each $i$. 

In this paper, when we use the term, surface-$(F_1,...,F_\mu)$-link, 
we suppose that each $F_*$ is a connected closed oriented surface unless otherwise specified. 
When we use the term, a surface-link $(L_1,...,L_\mu)$ 
or say that $(L_1,...,L_\mu)$ is a surface-link,  
 we suppose that  
each $L_*$ is a connected closed oriented surface  unless otherwise specified.

If all $F_i$ are spheres, $L$ is called a 
{\it 2-dimensional spherical link}, or {\it spherical link}. 
If $\mu=1$, $L$ is called a {\it surface-$F_1$-knot}, or {\it $F_1$-knot}.   
({\it 2-links} mean spherical links in some cases 
and 
surface-links in the other cases.)  
A {\it Seifert hypersurface} for $L$ is a connected compact oriented 3-dimensional submanifold $V\subset S^4$  such that $\partial V=L$. 
%
%

\begin{defn}\label{Idaho}
Surface-$(F_1,...,F_\mu)$-links,  
$L_0=(K_{01},...,K_{0\mu})$ and  $L_1=(K_{11},...,K_{1\mu})$,  
are said to be {\it surface-link-cobordant, link-cobordant}, or   {\it cobordant}, 
if there is an embedding map 
$$f:(F_1\x[0,1])\amalg...\amalg(F_\mu\x[0,1])\hookrightarrow S^4\x[0,1]$$ 
\noindent
with the following properties. 
Let $j\in\{0,1\}$. 
Let $\mu\in\N$. 
Let $i\in\{1,...,\mu\}$. 

\smallbreak\noindent (1)  
$f$ meets $S^4\x\{j\}$ transversely.  

$f^{-1}(S^4\x\{j\})$
$=(F_1\x\{j\})\amalg...\amalg(F_\mu\x\{j\}).$ 
 
\smallbreak\noindent (2)  
%
%
$f((F_1\x\{j\})\amalg...\amalg(F_\mu\x\{j\}))$ 
(resp. $f(F_i\x\{j\})$) in $S^4\x\{j\}$ is  $L_j$ (resp. $K_{ji}$).

\smallbreak
We say that an embedding map $f$ 
(resp. a submanifold 
$f((F_1\x[0,1])\amalg...\amalg(F_\mu\x[0,1]))\newline
\subset S^4\x[0,1]$) 
{\it gives cobordism} between $L_0$ and $L_1$.  
\end{defn}

A {\it genus $g$ handlebody} is 
a 3-dimensional compact connected oriented 3-manifold 
which consists of one 0-handle and $g$ copies of 1-handle. 
Let $\mu\in\N$. 
Let $i\in\{1,...,\mu\}$. 
Let $g_i$ be the genus of a closed oriented connected surface $F_i$. 
A surface-$(F_1,...,F_\mu)$-link  $L=(K_1,...,K_\mu)$ is called the {\it standard link} 
if there is a disjoint embedded 3-dimensional submanifold $(V_1,...,V_\mu)$ $\subset S^4$ 
such that  
$V_i$ is a genus $g_i$ handlebody
and such that  for each $i\in\{1,...,\mu\}$, 
$V_i$ is a Seifert hypersurface for $K_i$. 
%
%
A surface-$(F_1,...,F_\mu)$-link $L$ is said to be {\it trivial} if $L$ is spherical and standard.  
If a spherical link $L$ is cobordant to the trivial link, $L$ is said to be {\it slice}.

A surface-$(F_1,...,F_n)$-link  $L=(K_1,...,K_n)$ is called a {\it boundary link} 
if there is a disjoint embedded 3-dimensional submanifold 
$V_1\amalg...\amalg V_n\subset S^4$, where $\amalg$ denotes the disjoint union,  
such that  for each $i\in\{1,...,n\}$, 
$V_i$ is a Seifert hypersurface for $K_i$. 

Let  $L=(K_1,...,K_n)$  be a surface-$(F_1,...,F_n)$-link $\subset S^4$. 
Remove $K_{j_1},...,K_{j_p}$ from $L$. 
Suppose that the left components are $K_{l_1},...,K_{l_{n-p}}$. 
We do not suppose that \newline$l_1\leqq...\leqq l_{n-p}$.
We call the surface-link $(K_{l_1},...,K_{l_{n-p}})$ a {\it sublink} of $L$. 
We do not say that the empty set (resp. $L$ itself) is a sublink of $L$ in this paper.

\begin{prob}\label{outstanding}
Are all 2-dimensional spherical links slice?
\end{prob}

This is a well-known outstanding open problem. 
See \cite{CochranOrr, LevineOrr, Ogasadiscs, OgasaIntro,  Orr} 
for the history and the background.  
In order to attack this open problem in the future, 
in this paper we consider the following problem which includes the above one.

\begin{prob}\label{koremo}
Classify surface-links up to surface-link-cobordism.   
\end{prob}  

Why we consider Problem \ref{koremo} together with Problem \ref{outstanding}
is  
because we know many non-cobordant pairs in the case of all non-spherical surface-links 
(\cite{Cochran, OgasaZ, Orr, Ruberman, Sato}). 
We will review an outline of the results in the following paragraphs. 
We hope the following: 
we shall continue to make new pairs of non-cobordant, non-spherical  surface-links, and then we may solve Problem \ref{outstanding} in the future.  
Of course it itself is important to make such pairs. 
Our results,  Main Theorem \ref{Alabama} 
and Theorem \ref{main},   
show new examples which give answers to Problems \ref{koremo}.


Levine (unpublished), and Sato in \cite{Sato} defined   
the Sato-Levine invariant which is a surface-link-cobordism invariant  
and which is trivial for the standard link. 
\cite{Ruberman, Sato} showed a $(S^2,T^2)$-link 
whose Sato-Levine invariant is nontrivial. 
See \S\ref{South Carolina}.

In \cite{Cochran} Cochran defined Cochran sequence 
which is a surface-link-cobordism invariant 
and which is trivial for the standard link. 
He also proved that Cochran sequence is nontrivial for a $(S^2,T^2)$-link.  
See \S\ref{Vermont}.

In \cite{Orr} Orr defined Orr invariant 
which is a surface-link-cobordism invariant under a condition 
and which is trivial for the standard link. 
\cite[\S5]{Orr} claimed that Orr invariant is nontrivial for a surface-link.   
See \S\ref{Wyoming}. 
In \S\ref{Rhode Island} we prove a new property of Orr invariant (Theorem \ref{Virginia}). 


\cite[\S2]{Sato} proved that 
if two surface-links are cobordant and 
the alinking number of one of the two is zero, 
 then that of the other is zero. 
In \cite[Proposition 7.10]{OgasaZ} the author generalized it and proved that 
if two surface-links are cobordant 
 then the alinking number of the two are the same.  Note that 
the former result does not imply the latter one directly.  
%
It is well-known that the alinking number is nontrivial for a $(S^2,T^2)$-link 
and that it is trivial for the standard link. 
See \S\ref{Tennessee}.


%
%
%
In \cite[5.2 Theorem]{Stallings} Stallings proved the following:  
For a surface-link $A$,  
let $\pi_A$ be \newline 
$\pi_1(S^4-A)\cong\pi_1(S^4-N(A))\cong\pi_1(S^4-{\rm Int}N(A))$, 
where  
$N(A)$ is the tubular neighborhood of $A$ in $S^4$ and Int$N(A)$ the interior of $N(A)$.  
Let $\N_{\geqq2}=\N-\{1\}$. 
Let $k\in\N_{\geqq2}$. 
Let $G$ be a group.  
Let $G_k=[G,G_{k-1}]$ and $G_1=G$.  
If $G=\pi_L$, let $\pi_{L,k}$ denote $G_k$ (we use `,' in order to avoid confusion). 
Let $L$ and $L'$ be cobordant surface-links. 
It holds that  for any $k\in\N$, 
$\pi_L/\pi_{L,k} \cong \pi_{L'}/\pi_{L',k}$. 
In this paper we say that $\pi_L$ and $\pi_{L'}$ are {\it Stallings-equivalent}  
if  for any $k\in\N$, $\pi_L/\pi_{L,k}\cong\pi_{L'}/\pi_{L',k}$. 


%
%


%
%
%
In \S\ref{daidai} we submit a new way to detect surface-link-cobordism. 
We introduce terminologies, 
`$n$-covering-link ($n\in\N$)', 
`alinking-equivalent', 
and 
`weakly alinking-equivalent', 
We prove the following:   
Suppose that surface-links $\mathcal L$ and $\mathcal L'$ with a condition 
are cobordant and 
that $p$ is a prime power. Then 
 `the $p$-covering-link of $\mathcal L$' and  
 that of $\mathcal L'$  under a condition 
are  
`weakly alinking-equivalent'. 
See \S\ref{daidai} for the precise statement. 
This result is a theme of this paper. 

\bigbreak
We state our main result. 
There are two surface-links with the following properties. 
We cannot detect whether $\mathcal L$ and $\mathcal L'$ are non-cobordant  
by using Orr invariant, Cochran sequence, Sato-Levine invariant, 
the alinking number, 
and Stallings-equivalence which is defined as above.  
However our new way can detect it. 
See Main Theorem \ref{Alabama} for detail.

\begin{mth}\label{Alabama}  
There are $(S^2, S^2, T^2)$-links, $\mathcal L=(P,Q,R)$ and $\mathcal L'=(P',Q',R')$, 
with the following properties.  
Call all 2-component-sublinks of $\mathcal L$ $($resp.  $\mathcal L'),$ 
$\mathcal L_1=(Q, R)$,   
$\mathcal L_2=(P,R)$,  
$\mathcal L_3=(P,Q)$,  
$\mathcal L_{-1}=(R,Q)$,  
$\mathcal L_{-2}=(R,P)$,  and
$\mathcal L_{-3}=(Q,P)$ 
$($resp. 
$\mathcal L'_1=(Q,'R')$,   
$\mathcal L'_2=(P',R')$,  
$\mathcal L'_3=(P',Q')$,  
$\mathcal L'_{-1}=(R',Q')$,  
$\mathcal L'_{-2}=(R',P')$,  and
$\mathcal L'_{-3}=(Q',P').)$
Let $\mathcal S=\{1,2,3,-1,-2,-3\}$.

\begin{enumerate}
%
\item
Let $*\in\mathcal S$. 
The alinking number of $\mathcal L_*$
is the same as 
that of $\mathcal L'_*$.  

\item$({\rm i})$
The 3-component-Sato-Levine invariant of $\mathcal L$ is the same as 
that of $\mathcal L'$. 

\item[]$({\rm ii})$
Let $*\in\mathcal S$.  
The 2-component-Sato-Levine invariant of $\mathcal L_*$ is 
 the same as that of $\mathcal L'_*$.

\item
Let $*\in\mathcal S$. 
The Cochran sequence of $\mathcal L_*$ 
is the same as that of $\mathcal L'_*$. 

\item$({\rm i})$
%
$\pi_\mathcal L$ and $\pi_{\mathcal L'}$ are Stallings-equivalent. 
%

\item[]$({\rm ii})$
Let $*\in\mathcal S$. 
$\pi_{\mathcal L_*}$ and $\pi_{\mathcal L'_*}$ are Stallings-equivalent. 

\item$({\rm i})$
Let $\#\in\N_{\geqq2}\cup\{\omega\}$.   
The following two are equivalent. 

{\rm(I)} We can define the Orr invariant 
$\theta_\#(\mathcal L, \tau)$ for a meridian $\tau$ and for the element $\#$.  

{\rm(II)} We can define $\theta_\#(\mathcal L', \tau')$ for a meridian $\tau'$ and for the element $\#$.

Suppose that {\rm (I)} $($resp. {\rm (II)}$)$ holds. 
Then we have  $\theta_\#(\mathcal L)=\theta_\#(\mathcal L')=0$. 

\item[]$({\rm ii})$
Let $*\in\mathcal S$. 
We can define the Orr invariant 
$\theta_\#(\mathcal L_*, \alpha)$  
$($resp. $\theta_\#(\mathcal L'_*, \alpha')$$)$ 
for all meridians $\alpha$ 
$($resp.  $\alpha')$ 
and for all elements $\#\in\N_{\geqq2}\cup\{\omega\}$, 
and
we have $\theta_\#(\mathcal L_*)=\theta_\#(\mathcal L'_*)=0.$    







\item
Our new way in \S\ref{daidai} implies that $\mathcal L$ and $\mathcal L'$ are non-cobordant.  
\end{enumerate}
    
\end{mth}

\noindent{\bf Note.} 
We will define many links in this paper and we need many letters to represent them. 
So we use not only Roman-type letters but also  calligraphic letters for them.

\bigbreak
\section{The alinking number}\label{Tennessee}
\noindent   
The alinking number was introduced in \cite{Sato}. 
Let $l\in\N_{\geqq2}.$
Let $L=(K_1,...,K_l)$ be a surface-$(F_1,...,F_l)$-link in $S^4$. 
For any $i\in\{1,...,l\}$, 
take any circle embedded in $K_i$. Give any orientation to the circle.  
For any distinct $i,j\in\{1,...,l\}$, 
consider the linking number of the circle and $K_j$ in $S^4$. 
Make a set of all of the linking number. 
Then the set is regarded as $n\cdot\Z$ for a number $n\in\{0\}\cup\N$.  
Note that if $n=0$ then the set is $\{0\}$. 
We call this number $n$ the {\it alinking number} 
alk$(K_i\subset L, K_j\subset L)$  
of $K_i$ in $L$ around $K_j$ in $L$. 
Note that alk$(K_1\subset L, K_2\subset L)$ is not equal to 
alk$(K_2\subset L, K_1\subset L)$ in general. 

Let $L=(K_1,K_2)$ be a surface-$(S^2,F)$-link in $S^4$. 
Then alk$(K_1\subset L, K_2\subset L)=0$. 
Thus 
we let the alinking number of $L$ be alk$(K_2\subset L, K_1\subset L)$.

A surface-$(F_1,...,F_\mu)$-link  $L=(K_1,...,K_\mu)$ is called a {\it semi-boundary link} 
if  for each $i$, $[K_i]=0\in H_2(S^4-\amalg_{j\neq i, 1\leqq j\leqq\mu}K_j;\Z)$. 
For any distinct $i,j$, a Seifert hypersurface $V_i$ for $K_i$   
is called a {\it special Seifert hypersurface}   
if  $V_i\cap K_j=\phi$.  
It is well-known that 
a surface-$(F_1,...,F_\mu)$-link  $L=(K_1,...,K_\mu)$ is a semi-boundary link
if and only if 
there is a special Seifert hypersurface $V_i$ for $K_i$ for each $i$.  

\begin{pr}\label{South Dakota} {\bf {\rm(\cite[\S2]{Sato}.)}} 
%
%
Let $L=(K_1,K_2)$ be a surface-$(F_1,F_2)$-link. 
Then $(1) - (3)$ are equivalent each other. 
\smallbreak\noindent  
$(1)$ 
 $L$ is semi-boundary.   
\smallbreak\noindent  
$(2)$  
 $L$ is cobordant to a semi-boundary link. 
\smallbreak\noindent
$(3)$ 
Both ${\rm alk}(K_1\subset L,K_2\subset L)$ and 
${\rm alk}(K_2\subset L, K_1\subset L)$ are zero. 
\end{pr}

\begin{cor}\label{coro}{\rm(\cite[\S2]{Sato}.)}
If a surface-$(F_1,F_2)$-link $L=(K_1,K_2)$ is cobordant to the standard link,  
$L$ is semi-boundary.  
\end{cor}

Proposition \ref{South Dakota} implies that 
if two surface-links are cobordant and 
the alinking number of one of the two is zero, 
 then that of the other is zero. 
In \cite[Proposition 7.10]{OgasaZ} the author generalized it and proved the following. 
The proof is short so it is cited here.

\begin{thm}\label{alkcob}
If surface-$(F_a, F_b)$-links,  
$L_0=(K_{0a},K_{0b})$ and $L_1=(K_{1a},K_{1b})$,  are cobordant, 
then 
$${\rm alk}(K_{0a}\subset L_0, K_{0b}\subset L_0)=
{\rm alk}(K_{1a}\subset L_1, K_{1b}\subset L_1)$$
and 
$${\rm alk}(K_{0b}\subset L_0, K_{0a}\subset L_0)=
{\rm alk}(K_{1b}\subset L_1, K_{1a}\subset L_1).$$
\end{thm}

\noindent
{\bf Proof of Theorem \ref{alkcob}.} 
Take $S^4\x[0,1]$. 
Let  $j\in\{0,1\}$. 
Regard $L_j$ as a submanifold of $S^4\x\{j\}$. 
Take an embedding map 
$f:(F_a\x[0,1])\amalg(F_b\x[0,1])\hookrightarrow S^4\x[0,1]$ 
which gives cobordism between $L_0$ and $L_1$.  
Let $V_{jb}$ be a Seifert hypersurface for $K_{jb}\subset S^4\x\{j\}$, where 
 $K_{ja}\cap V_{jb}\neq\phi$ may hold.  
Note that there is an embedded compact oriented 4-manifold $W_b\subset S^4\x[0,1]$ 
whose boundary is 
$V_{0b}\cup f(F_b\x[0,1])\cup V_{1b}$. 
Take any embedded circle  $C\subset F_a$ and any $t\in[0,1].$  
By using $f(C\x\{t\})\cap W_b$, 
we can prove that \newline
${\rm alk}(K_{0a}\subset L_0, K_{0b}\subset L_0)=
{\rm alk}(K_{1a}\subset L_1, K_{1b}\subset L_1)$
holds. 
By replacing $a$ with $b$ in the above proof, 
we can prove that 
${\rm alk}(K_{0b}\subset L_0, K_{0a}\subset L_0)=
{\rm alk}(K_{1b}\subset L_1, K_{1a}\subset L_1)$ holds. 
\qed

\bigbreak
\section{The Sato-Levine invariant}\label{South Carolina}
\noindent  
The Sato-Levine invariant was defined 
by Sato (\cite{Sato}) and by Levine (unpublished). 
Let $L=(K_1,K_2)$ be a semi-boundary surface-$(F_1,F_2)$-link. 
For each $i\in\{1,2\}$, 
let $V_i$ be a special Seifert hypersurface for $K_i$. 
We can suppose 
that $V_1$ intersects $V_2$ transversely. 
Note that $G=V_1\cap V_2$ is a (not necessarily connected) closed oriented surface. 
By using $G$, $V_1$ and $V_2$,  we can define Thom-Pontrjagin map
$p:S^4\to S^2$ such that the inverse image of a regular value is $G$.   
The {\it $($2-component$)$ Sato-Levine invariant} $\beta(L)$ is $[p]\in\pi_4(S^2)\cong\Z_2$. 
In \cite{OgasaSL} the author proved that we can regard $\beta(L)$ as an element 
$\in\Omega^{\rm{spin}}_2\cong\Z_2$, 
where $\Omega^{\rm{spin}}_2$ is the second spin cobordism group.  
We can define the {\it 3-component Sato-Levine invariant} $\beta(L)\in\pi_4(S^3)\cong\Z_2$ for any semi-boundary surface-$(F_1, F_2, F_3)$-link $L=(K_1,K_2, K_3)$.   
We can define the {\it 4-component Sato-Levine invariant} $\beta(L)\in\pi_4(S^4)\cong\Z$ for any semi-boundary surface-$(F_1, F_2, F_3, F_4)$-link $L=(K_1,K_2, K_3, K_4)$.   
If $m\geqq5$, it is nonsense to define the Sato-Levine invariant for $m$-component surface-link $L$. 

\cite{Ruberman, Sato} proved that 
there is a $(S^2,T^2)$-link whose Sato-Levine invariant is nontrivial.


In \cite{Orr} Orr proved outstanding results: 
Let  $m\in\{2,3,4\}$. If $L$ is an $m$-component  spherical 2-link, the Sato-Levine invariant of $L$ is zero (\cite[Theorem 4.1]{Orr}).   
If $L$ is a 4-component (not necessarily spherical) surface-link, 
the Sato-Levine invariant of $L$ is zero  (\cite[Theorem 3.5]{Orr}).



\bigbreak
\section{The Cochran sequence}\label{Vermont}
\noindent
The Cochran sequence was defined in \cite{Cochran}. 
Let $*\in\N$. Let $F$ (resp. $G$, $F_*$)  
be a connected closed oriented surface through this section. 

\smallbreak
\noindent{\bf Definition.}  
A $(S^2, F)$-link $L=(K_1, K_2)\subset S^4$ and 
a $(S^2, G)$-link $L'=(K'_1, K'_2)\subset S^4$ 
are {\it weak-cobordant} if they satisfy the following conditions:
 
\smallbreak\noindent(1)
The alinking number of $L$ (resp. $L'$) is zero. 
Hence there are Seifert hypersurfaces $Z$ (resp. $Z'$) for $K_1$ (resp. $K'_1$) 
in $S^4-{\rm Int}N(K_2)$ (resp.  $S^4-{\rm Int}N(K'_2)$), 
where $N(K_*)$ is the tubular neighborhood of $K_*\subset S^4$. 

\smallbreak\noindent(2)   
Let $I=[0,1]$. 
There is a compact oriented 3-dimensional submanifold $W\amalg Y\subset S^4\x I$
such that 

\smallbreak\noindent\hskip1cm(a) 
$\partial(W\amalg Y)=L\amalg L'$ 

\smallbreak\noindent\hskip1cm(b) 
$W=S^2\x I$

\smallbreak\noindent\hskip1cm(c) 
The closed compact oriented 3-dimensional submanifold $Z\cup W\cup Z'$ 
bounds a

\noindent\hskip15mm
compact oriented 
4-dimensional submanifold $Q$ $\subset\overline{(S^4\x I)-N(Y)}$, 

\noindent\hskip15mm
where $N(Y)$ is the tubular neighborhood of $Y\subset S^4\x I$. 

\noindent\hskip15mm
$Q\cap(S^4-{\rm Int}N(K_2))=Z$. 
$Q\cap(S^4-{\rm Int}N(K'_2))=Z'$. 

\bigbreak
It is trivial that the following holds.

\begin{pr}\label{Utah}
If $(S^2, T^2)$-links,  $L=(K_1, K_2)$ and $L'=(K'_1, K'_2)$, $\subset S^4$  
are cobordant 
and if the alinking number of $L$ is zero, 
 then they  are weak-cobordant. 
\end{pr}

\noindent{\bf Definition.}  
Let $L=(K_1, K_2)$ be a $(S^2, T^2)$-link whose alinking number is zero. 
By \cite[Theorem 4.1]{Cochran} there is a Seifert hypersurface $V_i$ for $K_i$ ($i=1,2$) 
such that $V_1$ intersects $V_2$ transversely 
and such that $V_1\cap V_2=J$ is a closed oriented connected surface.   
The {\it derivative link} $D(L)$ of $L$ is 
the weak-cobordism class of the surface-link $(K_1, J)$. 

\bigbreak\noindent
{\bf Theorem. (\cite[Theorem 4.2.]{Cochran}.)}
{\it 
If  
a $(S^2, F)$-link $L=(K_1, K_2)\subset S^4$ and 
a $(S^2, G)$-link $L'=(K'_1, K'_2)\subset S^4$  
are weak-cobordant, 
then 
$D(L)$ and $D(L')$ are weak-cobordant.  
} 

\bigbreak\noindent
{\bf Theorem. (\cite[Corollary 4.3.]{Cochran}.)}
{\it 
For any $(S^2, F)$-link $L=(K_1, K_2)$ whose alinking number is zero, 
the sequence of links 
$$D^0(L)=L, D(L), D^2(L)=D(D(L)), D^3(L)=D(D(D(L)))...$$ 
is well defined in the category of links modulo weak-cobordant.  
} 

\bigbreak\noindent{\bf Definition.} 
Let $L=(K_1, K_2)$ be a $(S^2, F)$-link whose alinking number is zero. 
Note that for any natural number $n$,  
$D^n(L)$ is a $(S^2, F_n)$-link whose alinking number is zero. 
Note that we can define the Sato-Levine invariant 
$\beta(D^n(L)) (n\in\N)$. 
Let $\beta^{n+1}(L)=\beta(D^n(L)) (n\in\N)$ and $\beta^1(L)=\beta(L)$. 
The {\it Cochran sequence} for $L$ is 
the sequence $\{\beta^n(L)\}_{n\in\N}$. 

\bigbreak\noindent
{\bf Theorem. (\cite[Corollary 5.2]{Cochran}.)}
{\it 
Let $L=(K_1, K_2)$ be a $(S^2, F)$-link whose alinking number is zero. 
The Cochran sequence $\{\beta^n(L)\}_{n\in\N}$ 
is weak-cobordism invariant. 
}

\bigbreak
It is trivial that the following holds.

\begin{thm}\label{Texas}
Suppose that $(S^2, T^2)$-links, $L$ and $L'$, are cobordant 
and that the alinking number of $L$ is zero. 
Then the Cochran sequences, $\{\beta^n(L)\}_{n\in\N}$ and $\{\beta^n(L')\}_{n\in\N}$, 
are equivalent.  
\end{thm}

\cite{Cochran} proved that there is a $(S^2, T^2)$-link with a nontrivial Cochran sequence.






\bigbreak
\section{The Orr invariant}\label{Wyoming}
\noindent  
%
%
In \cite[\S2]{Orr} Orr defined 
the Orr invariant $\theta_k$ for any $k\in\N_{\geqq2}$ and $\theta_\omega$ 
in Definition \ref{Wisconsin} 
for codimension two closed oriented submanifolds $L\subset S^{n+2}$     
if $L$ is a disjoint union of connected homology spheres.   
He used his invariant and proved \cite[Theorem 4.1]{Orr},  
which is explained in \S\ref{South Carolina}.  
Furthermore, in \cite[the first few lines of \S5]{Orr},  
he stated that 
$\theta_k$ and $\theta_\omega$ can be defined in some other cases 
which include ones in Definition \ref{Wisconsin}. 


\begin{defn}\label{Wisconsin}  
Let $m\in\N_{\geqq2}$. 
Let $L=(K_1,...,K_m)$ be an $m$-component surface-link. 
Let $E_L=S^4-\text{Int}N(L)$. 
Note that $N(L)=L\x D^2$. 
Note that $\partial E_L=L\x S^1=\partial N(L)$.

Let $p$ be a base point of $S^4$. 
Let $S_1^1\vee...\vee S_m^1$ be a bouquet. 
Let $b=S_1^1\cap...\cap S_m^1$ be a base point of $S^1\vee...\vee S^1$. 
Let $(L,\tau)$ be an $m$-component surface-link with a fixed meridian. 
The meridian of $L$ defines a continuous map 
$$\tau:(S_1^1\vee...\vee S_m^1, b)\to(E_L, p).$$ 
Here, we use $\tau$ again.  
Let $i\in\{1,...,m\}$. 
Note that $\tau|_{S_i^1}$ defines the meridian of $K_i$.    

Let $\pi_L$ denote $\pi_1E_L$. 
Let $F=\underbrace{\Z*...*\Z}_{m}$. 
%
%
Suppose that the condition (1) (resp. (2)) holds. 

\smallbreak\noindent(1)  
$\pi_L/\pi_{L.k}=F/F_k$ for a natural number $k$. 
 
\smallbreak\noindent(2)  
%
%
%
$\displaystyle\lim_{\infty \leftarrow k}\pi_L/\pi_{L.k}$ 
$=$
$\displaystyle\lim_{\infty \leftarrow k}F/F_k$ 
%
%

\smallbreak
\noindent
(By \cite[Hopf's theorem]{Brown} and \cite{Stallings}, 
we have the following:  if $L$ is a disjoint union of connected homology spheres, 
the congruence in (1) holds for all natural numbers $k$  
and the condition (2) 
holds.)  

\smallbreak

Let $K(G,1)$ be the Eilenberg-MacLane space for any group $G$. 
Let $\bar\tau$ denote the induced homotopy-type-equivalence map 
$K(F/F_k,1)\to K(\pi_L/\pi_{L,k},1)$  
(resp.  
$K(F/F_\omega,1)\newline \to K(\pi_L/\pi_{L,\omega},1)$). 
Let $\phi_k:E_L\to K(F/F_k,1)$  
(resp. $\phi_\omega:E_L\to K(\displaystyle\lim_{\infty \leftarrow k}F/F_k,1)$) 
be the composition of 
`the map 
$E_L\to K(\pi_L/\pi_{L,k},1)$ (resp. $E_L\to K(\pi_L/\pi_{L,\omega},1)$)
realizing 
$\pi_L\to\pi_L/\pi_{L,k}$ (resp. $\pi_L\to\pi_L/\pi_{L,\omega}$)'  
with the map $(\bar\tau)^{-1}$ 
Thus we define the following homomorphism naturally: 
$$\pi_1K_i\to\pi_1K_i\x\pi_1S^1\cong\pi_1(K_i\x S^1)\to\pi_L\to\pi_L/\pi_{L,k}\to F/F_k$$
$$\text{(resp.}    
\pi_1K_i\to\pi_1K_i\x\pi_1S^1\cong\pi_1(K_i\x S^1)\to
\pi_L\to
\displaystyle\lim_{\infty \leftarrow k}\pi_L/\pi_{L,k} 
\to
\displaystyle\lim_{\infty \leftarrow k}F/F_k
\text{)}$$
This induces a homomorphism 
$\pi_1K_i\to F/F_k$ 
(resp. $\pi_1K_i\to\displaystyle\lim_{\infty \leftarrow k}F/F_k$).   
Suppose that we have the condition 

\smallbreak\noindent(3)  
$\pi_1K_i\to F/F_k$ 
(resp. 
$\pi_1K_i\to
\displaystyle\lim_{\infty \leftarrow k}F/F_k$)  
is the zero map. 

\smallbreak\noindent
(If $L$ is a disjoint union of connected homology spheres, 
 the condition (3) holds. 
See \cite[the part from the last line of page 546 to the first line of page 547]{Orr}. )

\smallbreak
The quotient homomorphism 
$$\psi_k:F\to F/F_k$$
induces an inclusion of Eilenberg-MacLane spaces
$$\psi_k:S^1_1\vee...\vee S^1_m=K(F,1)\to K(F/F_k,1).$$
Let $K_k$ be the mapping cone of $\psi_k$. Note that $K_k$ is simply-connected.

The homomorphisms $\psi_k$  induce a homomorphism 
$$\psi_\omega:F\to
\displaystyle\lim_{\infty \leftarrow k}F/F_k
$$
and an inclusion of the spaces 
$$\psi_\omega:S^1_1\vee...\vee S^1_m=K(F,1)\to 
K(\displaystyle\lim_{\infty \leftarrow k}F/F_k, 1).$$
Let $K_\omega$ be the mapping cone of $\psi_\omega$.

%
The condition (3) 
and the property of the Eilenberg-MacLane space 
imply the following commutative diagram.

$$\begin{matrix}
K_i\x S^1  & \stackrel{\phi_k|_{K_i\x S^1}}\to & K(F/F_k,1) \\
                                          
\downarrow_{\text{projection}}  &\hskip-10mm\circlearrowright 
& & & \\
                                          &\nearrow & \\
S^1 &  & \\
\end{matrix}$$


Therefore
$\phi_k$ (resp. $\phi_\omega$) extends canonically to 
a continuous map 
$$\bar{\phi_k}:S^4\to K_k$$   
$$({\rm resp.} \bar{\phi_\omega}:S^4\to K_\omega).$$


\noindent Define the {\it Orr invariant} $\theta_k(L,\tau)$ (resp. $\theta_\omega(L,\tau)$) 
to be 
$[\bar{\phi_k}]\in \pi^4(K_k)$ (resp. $[\bar{\phi_\omega}]\newline\in\pi^4(K_\omega)$),  
where $k\in\N_{\geqq2}$.

If $\theta_*$ is defined for $*\in\Lambda$ which is an infinite set, 
the sequence $\{\theta_*\}_{*\in\Lambda}$ is called {\it Orr sequence}.  
\end{defn}

\cite[Theorem in \S5]{Orr} claimed 
that there is a surface-link with a nontrivial Orr invariant. 

We have the following theorem.  
\begin{thm}\label{Atlanta}
{\bf(\cite[Theorem 2.1.(iv)]{Orr}.)} 
Let $k\in\N_{\geqq2}\cup\{\omega\}$. 
We have the following: 
$\theta_k(L,\tau)=0$ for a choice $\tau$ of meridians for $L$ 
if and only if 
$\theta_k(L,\tau')=0$ for any other choice $\tau'$ of meridians for $L$.  
\end{thm}

By Theorem \ref{Atlanta}, the notation, `$\theta_k(L)=0$', makes sense. 
(Note that we omit a meridian in $(\quad)$ of $\theta_k(\quad)$.)

\section{Boundary links and spun surface-links}\label{Rhode Island}

\noindent 
We use boundary links and spun surface-links in order to make examples in the proof. 
We defined boundary links in \S\ref{Introduction}. 
We now define spun surface-links. 

\begin{defn}\label{Georgia}
We say that a surface-link $L$ is a {\it spun surface-link} 
if $L$ is (not necessarily order-preserving) equivalent to  
the surface link $L'$ which is obtained as follows. 

Regard $S^4$ as $\R^4\cup{\infty}$. 
Regard $\R^4=\{(x,y,z,w)|x, y, z,w\in\R\}$. 
Regard 
$\R^4\newline=\{(x,y,z,w)|x, y, z,w\in\R\}$ 
as 
the result of rotating  
$Q\newline=\{(x,y,z,w)|x, y\in\R, z\geqq0,w=0\}$ 
around 
$A=\{(x,y,z,w)|x, y\in\R, z=0,w=0\}$.

Let $I_1,...,I_p$ be the intervals. 
Let $f:I_1\amalg...\amalg I_p\amalg S^1_1\amalg...\amalg S^1_q\hookrightarrow Q$ 
be an embedding map 
such that for each $*$, 
$f(I_*)\cap A$ is the two points $f(\partial I_*)$.   
Rotate Im$f$ around $A$ together when we rotate $F$ around $A$ and obtain $\R^4$. 
Thus we obtain a surface-link $L'\subset S^4$.   
\end{defn} 

We prove some properties of boundary links and spun surface-links.

\begin{pr}\label{standard}
Let $L=(K_1,...,K_n)$ be a boundary surface-$(F_1,...,F_n)$-link. 
Then $L=(K_1,...,K_n)$ is cobordant to the standard link.  
\end{pr}

\noindent
{\bf Proof of Proposition \ref{standard}.} 
Let $L'=(K'_1,...,K'_n)$ be the standard surface-$(F_1,...,F_n)$-link.   
Embed a 4-ball $B$ (resp. $B'$)  in  $S^4$ such that $B\cap B'=\phi$. 
Take $L$ (resp. $L'$) in $B$ (resp. $B'$).   
Let $V_1\amalg...\amalg V_n\subset B$ (resp. $V'_1\amalg...\amalg V'_n\subset B'$) 
be a disjoint embedded 3-manifold 
 such that 
$V_i\cap V_j=\phi$ (resp. $V'_i\cap V'_j=\phi$)   
for 
all distinct $i, j$,  
and such that 
$V_i$ (resp. $V'_i$)
is a Seifert hypersurface for 
$K_i$ (resp. $K'_i$)
for each $i$.    
Suppose that 
 $V'_i$ is an oriented genus $g_i$ handlebody for each $i$, where $g_i$ is the genus of $F_i$. 
Connect $V_i$ and $V'_i$ 
by an embedded 3-dimensional 1-handle $h^1_i\subset S^4$ for each $i$ 
such that $h^1_i\cap h^1_j=\phi$ 
for all distinct $i, j$, 
and call the result $W_i$. 
See \cite{Browder, Kirby, Wall} for handles, surgery and their associated terms. 
Note that the boundary of $W_1\amalg...\amalg W_n$ is a surface-$(F_1\sharp F_1,...,F_n\sharp F_n)$-link, 
where $\sharp$ denotes the connected-sum.

\begin{figure}
\includegraphics[width=11cm]{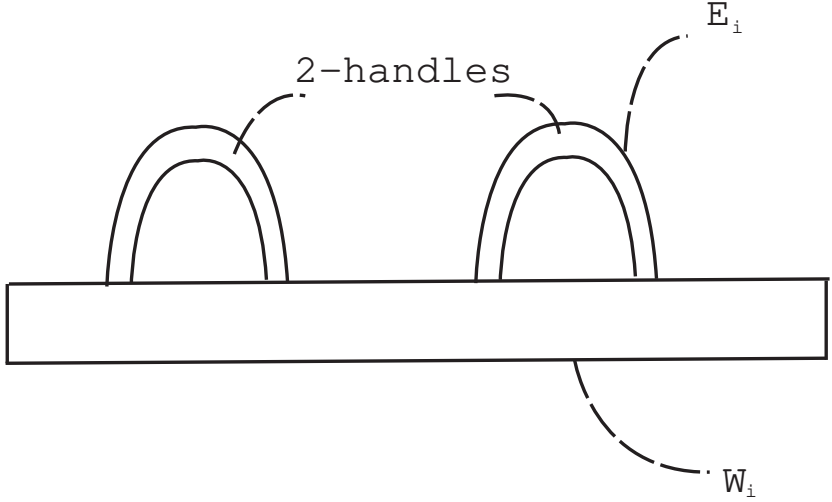}

\caption{\bf  An example of the handle decomposition of $O_i$ \label{Arizona}}

\end{figure}

Since $\Omega^{\rm Spin}_3\cong0$ (see \cite{Kirby}), 
there is a compact oriented spin 4-manifold $O_i$ 
with  a handle decomposition (see Figure \ref{Arizona}) 

\smallbreak
\hskip25mm
$O_i=(W_i\x[0,1])\cup$(4-dimensional 2-handles)$\cup(E_i\x[0,1])$     

\smallbreak
\noindent
to satisfy the following conditions: 
$W_i$ is the bottom and $E_i$ is the top. 
$\partial O_i\newline =(W_i\x\{0\})\cup(E_i\x\{1\})$. 
$E_i=(F_i-\text{(an embedded open 2-disc))}\times [0,1]$. 
     
Note that, here, we make $O_i$ as an absolute manifold not as a submanifold. 
We make a submanifold which is diffeomorphic to $O_i$ from now.


Recall that $L$ and $L'$ are in $S^4$. 
Take $S^4\x[0,1]$. Take $(S^4, L\amalg L')\x[0,1]$, 
where $(X,Y)$ denotes a pair of  a manifold $X$ and a submanifold $Y\subset X$. 
Take the tubular neighborhood $N(W_i)=W_i\x[-1,1]$ of $W_i$ in $S^4$. 
Take  
$(S^4, (W_1\x[-1,1])\amalg...\amalg(W_n\x[-1,1]))\x[0,1]$. 
Take $O_i\times[-1,1]$.  
Make these $S^4\x[0,1]$ and $O_i\times[-1,1]$ into a 5-manifold as follows:  
Identify $W_i\x[-1,1]\x[0,1]\subset S^4\x[0,1]$ with 
 $W_i\x[0,1]\times[-1,1] \subset O_i\times[-1,1]$
so that for each $s\in[-1,1]$ and each $t\in[0,1]$, 
$W_i\x\{s\}\times\{t\}\subset S^4\x[0,1]$ coincides with 
$W_i\x\{t\}\times\{s\}\subset O_i\times[-1,1]$. 
%
%
%
%
Thus the resultant 5-manifold $M$ 
has a handle decomposition 
\smallbreak 
\hskip25mm   
$M=(S^4\x[0,1])\cup\text{(5-dimensional 2-handles)}.$ 
\smallbreak 
\noindent 
Since $M$ is a spin oriented 5-manifold, 
\smallbreak 
\hskip25mm 
$M\cong(\natural^\mu S^2\x B^3)-\text{(an open 5-ball)}$,   
\smallbreak
\noindent 
where $\natural$ denotes the boundary-connected-sum and  $\mu$ is a natural number.

Note that $M$ can be embedded in $B^5$ so that $S^4\x\{0\}$ coincides with $\partial B^5$. 
Recall that 
$L\subset B\subset S^4=\partial B^5\subset B^5$
and that 
$L'\subset B'\subset S^4=\partial B^5\subset B^5$. 
Recall that $E_i\cup h^1_i\subset B^5$ 
and that  $E_i\cup h^1_i$ is diffeomorphic to $F_i\x[-1,1]$.   
Therefore 
we can embed $F_i\x[-1,1]\subset B^5$ 
such that 
$F_i\x\{-1\}$ 
(resp. $F_i\x\{1\}$) coincides with 
$K_i$ of $L\subset B\subset S^4$  (resp. $K'_i$ of $L'\subset B'\subset S^4$) for each $i$ 
and such that 
$(F_i\x[-1,1])\cap(F_j\x[-1,1])=\phi$ 
for all distinct $i, j$.
Hence Proposition \ref{standard} holds. \qed

\bigbreak
By Proposition \ref{standard}, 
Theorem \ref{chau}.(1) $\Rightarrow$ Theorem \ref{chau}.(2) is trivial.   
However the converse is not true in general.

\begin{thm}\label{chau}
Let $L=(K_1,...,K_n)\subset S^4$ be a surface-$(F_1,...,F_n)$-link. 
Then 
$(1)\Rightarrow(2)$ is true,  
but 
$(2)\Rightarrow(1)$ is false in general.  

\smallbreak
\noindent
$(1)$  
$L$ is cobordant to a boundary-link.

\smallbreak
\noindent
$(2)$  
Suppose that $L\subset S^4=\partial B^5\subset B^5$.  
Let $g_i$ be the genus of $F_i$. 
There is an embedded submnaifold 
$O_1\amalg...\amalg O_n\subset B^5$ with the following properties. 
Let $i,j\in\{1,...,n\}.$ 

\smallbreak\noindent$({\rm i})$
$O_i$ is an oriented genus $g_i$ handlebody. 

\smallbreak\noindent$({\rm ii})$
$O_i\cap O_j=\phi$ 
for all distinct $i, j$.

\smallbreak\noindent$({\rm iii})$
$O_1\amalg...\amalg O_n$ meets $\partial B^5$ transversely. 
$O_i\cap\partial B^5=\partial O_i$.  

\hskip2mm
$\partial O_1\amalg...\amalg\partial O_n$ $($resp. $\partial O_i)$ 
in $\partial B^5$ is 
$L$ $($resp. $K_i).$  
\end{thm}

\noindent
{\bf Note.}  
It is trivial that if $g_i=0$ for each $i$, $(2)\Rightarrow(1)$ is true.


\bigbreak\noindent
{\bf Proof of Theorem \ref{chau}.} 
We prove that Theorem \ref{chau}.(2) $\Rightarrow$ Theorem \ref{chau}.(1) is false in general.  
We define a $(S^2, T^2)$-link $L=(K, J)$ as follows: 
Let $K$ be the trivial spherical 2-knot $\subset S^4$. 
Let $N(K)$ be the tubular neighborhood of $K$ in $S^4$. 
Note that $N(K)=K\x D^2$. 
Let $C$ be an embedded circle $\{*\}\x\partial D^2\subset S^4$. 
Embed a solid torus $S^1\x D^2$ so that $S^1\x\{\sharp\}$ coincides with $C$. 
Let $J$ be an embedded torus $\partial(S^1\x D^2)$. 
%
By the construction, $L$ satisfies Theorem \ref{chau}.(2).  
By the construction, alk$(J,K)$=1. 
By Proposition \ref{South Dakota}  and Theorem \ref{alkcob},  
$L$ does not satisfy Theorem \ref{chau}.(1).  
 \qed

\begin{figure}
\includegraphics[width=13cm]{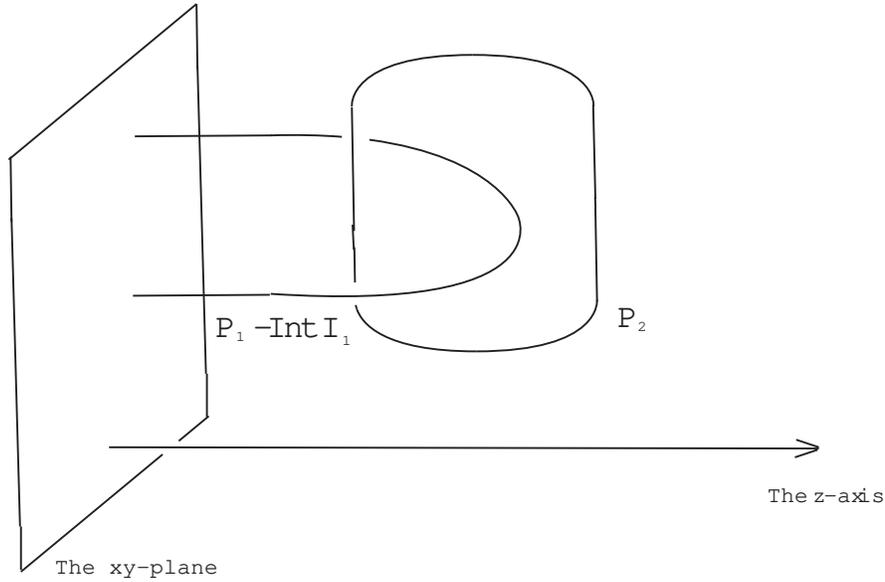}

\bigbreak
\caption{\bf  
A 1-link $(P_1, P_2)$ whose spun surface-link is  
a $(S^2, T^2)$-link with a nontrivial alinking-number.  
 \label{Arkansas}}

\end{figure}

\begin{note}\label{North Dakota}  
(1) 
Take a spun $(S^2, T^2)$-link $P$ of Figure \ref{Arkansas}.    
There, we omit the $x$-, $y$-axes and  draw the $z$-axis. 
%
%
We claim that  $P$ is equivalent to $L$ in Proof of Theorem \ref{chau}.   
{\it Reason.} 
Let ($S^4$ in Theorem \ref{chau})$-$(a point $\{\infty\}$) 
be $\R^4$ in Proof of Theorem \ref{chau}. 
Move $C$ in Proof of Theorem \ref{chau} in $S^4$ by an isotopy, and let $C$ include $\{\infty\}$. 
Then $L$  becomes $P$.

\smallbreak
\noindent 
(2)
Not all spun 
$(S^2, T^2)$-links are cobordant to the standard link  
because of the existence of 
$P$ in Note \ref{North Dakota}.(1).    
On the contrary it holds that all spun $m$-component spherical 2-links are slice, that is,  cobordant to the trivial link ($m\in\N$). 
Furthermore we will prove in \S\ref{New Jersey} that 
not all semi-boundary spun $(S^2,T^2)$-links are cobordant to the standard link. 
\end{note}

\begin{thm}\label{Virginia}
Let $L$ be a spun surface-link. 
If for a meridian $\tau$ and for an element $k\in\N_{\geqq2}\cup\{\omega\}$, 
we can define the Orr invariant $\theta_k(L, \tau)$,  
then we have $\theta_k(L)=0.$   
%
%
%
%
\end{thm}

\begin{lem}\label{Oklahoma}
Let $L$ be a spun surface-link $\subset S^4$.
Then there is an embedding map 
$g:D^3\amalg...\amalg D^3\amalg(D^2\x S^1)\amalg...\amalg(D^2\x S^1)
\hookrightarrow B^5$ with the following properties: 

\smallbreak\noindent$(1)$ 
$g$ meets $\partial B^5$ transversely. 
Call 
$D^3\amalg...\amalg D^3\amalg(D^2\x S^1)\amalg...\amalg(D^2\x S^1)$, $Y$. 
$g^{-1}(\partial B^5)$ is 
$\partial Y$.
$g(\partial Y)$ in $\partial B^5$ is 
a surface-link which is $($not necessarily order-preserving$)$ equivalent to $L$.

\smallbreak\noindent$(2)$  
Call the submanifold $g(Y)\subset B^5$,  $M$.  
Let $E_M=\overline{B^5-N(M)}$, 
where $N(M)$ is the tubular neighborhood of $M$ in $B^5$ 
and $\overline{{\color{white}abc}}$ denotes the closure in $B^5$.  
Let $\pi_M$ denote $\pi_1(E_M)$. 
Recall that $E_L=S^4-{\rm Int}N(L)$ and $\pi_L$ which are defined in \S\ref{Wyoming}. 
The inclusion map 
$E_L\hookrightarrow E_M$ 
induces an isomorphism map $\pi_L\cong\pi_M$. 
\end{lem}




%
%
%
%


\bigbreak\noindent{\bf Proof of Lemma \ref{Oklahoma}.}   
Take an embedding map 
$f:I_1\amalg...\amalg I_p\amalg S^1_1\amalg...\amalg S^1_q\newline
\hookrightarrow Q=\{(x,y,z,w)|x, y\in\R, z\geqq0,w=0\}$,  
which makes a surface-link 
which is $($not necessarily order-preserving$)$ equivalent to $L$,   
as in Definition \ref{Georgia}   for appropriate natural numbers $p$ and $q$. 
%
%
%
%
As we stated there, 
we regard  $\R^4=\{(x,y,z,w)|x, y, z,w\in\R\}$ 
as the result of rotating  
$Q=\{(x,y,z,w)|x, y\in\R, z\geqq0,w=0\}$ 
around 
$A\newline=\{(x,y,z,w)|x, y\in\R, z=0,w=0\}$.  
%

Let $\R^4_{z\geqq0}=\{(x,y,z,w)|z\geqq0\}$. 
Recall that $L\subset\R^4$.  
Let $L_{z\geqq0}=L\cap\R^4_{z\geqq0}$.  
Regard $\R^5_{v\geqq0}=\{(x,y,z,w,v)|v\geqq0\}$ 
as the result of rotating  $\R^4_{z\geqq0}$ through 180 degrees 
around $\{(x,y,z,w)|z=0\}$ as the axis. 
When we carry out this half rotation, we rotate $L_{z\geqq0}$  
together. 
Call the resultant submannifold contained in $\R^5_{v\geqq0}$,  $M$. 
Note that we have the following: 
$M$ is a disjoint union of $p$ copies of $D^3$ and $q$ copies of $D^2\x S^1$. 
$M\newline\cap\{(x,y,z,w,v)|v=0\}=\partial M$.  
$\partial M$ in $\{(x,y,z,w,v)|v=0\}\cup\{\infty\}$ is $L$ if we give an appropriate order to its components.

By the construction, 
$\overline{Q-N({\rm Im} f)}$ 
is homotopy type equivalent to 
$\overline{\R^4_{z\geqq0}-N(L_{z\geqq0})}$, 
and to 
$\overline{\R^5_{v\geqq0}-N(M)}$,  
and 
each homotopy-type-equivalence map is 
each natural inclusion. 
By using this homotopy type equivalence, Lemma \ref{Oklahoma} holds. \qed\bigbreak

\noindent{\bf Proof of Theorem \ref{Virginia}.} 
By the assumption,  we can define the Orr invariant $\theta_k(L, \tau)$  
by using a continuous map 
$\bar\phi_k:S^4\to K_k$ as in \S\ref{Wyoming}. 
Recall that the map $\bar\phi_k$ is made from 
$\phi_k:E_L\to K(F/F_k,1)$  (or $\displaystyle\lim_{\infty \leftarrow k}K(F/F_k,1)$).

Call $f(\partial Y)\subset \partial B^5$, $L$ again, and 
let $L=(L_1,...,L_r)$, where $r=p+q$. 
Let $*\in\{1,...,r\}$. 
Let $M$ be a disjoint union of compact connected components $M_1\amalg...\amalg M_r$ 
such that $\partial M_*=L_*$. 
Consider the following commutative diagram of continuous maps. 
The maps are inclusion maps

\smallbreak
\hskip2cm
$\begin{tabular}{ccccc}
$M_i\x\partial D^2$&$\to$&$M\x\partial D^2$&$\to$&$E_M$\\
$\uparrow$& &$\uparrow$& &$\uparrow$\\
$\partial N(L_i)=L_i\x\partial D^2$ & $\to$ & $\partial E_L$ & $\to$ & $E_L$\\
\end{tabular}$

\smallbreak\noindent
This induces the following commutative diagram of homomorphisms. 
The homomorphism represented by the left uparrow is called $\xi$. 

\smallbreak
\hskip2cm
$\begin{tabular}{ccccc}
$\pi_1(M_i\x\partial D^2)$ & $\to$ & $\pi_1(M\x\partial D^2)$ & $\to$ & $\pi_M$\\
$\uparrow_{\xi}$& &$\uparrow$& & $\uparrow_{\cong}$\\
$\pi_1(L_i\x\partial D^2)$ & $\to$ & $\pi_1(\partial E_L)$ & $\to$ & $\pi_L$\\
\end{tabular}$

By the assumption and Lemma \ref{Oklahoma}, 
the homomorphism of the right uparrow $\pi_M\newline\to\pi_L$ is isomorphic. 
Therefore the map $\phi_k$ extends to a map $\Phi_k$ 
$E_M\to K(F/F_k,1)$  (or $\displaystyle\lim_{\infty \leftarrow k}K(F/F_k,1)$).

We have $\pi_1(L_i\x\partial D^2)=\pi_1(L_i)\x\pi_1(\partial D^2)$
and 
$\pi_1(M_i\x\partial D^2)=\pi_1(M_i)\x\pi_1(\partial D^2)$. 
It holds that $\xi(\pi_1(L_i))\subset\pi_1(M_i)$. 
Furthermore 
this inclusion map $\pi_1(L_i)\to\pi_1(M_i)$ is onto. 
Therefore the condition (3) in \S\ref{Wyoming} 
and the property of the Eilenberg-MacLane space 
imply the following commutative diagram.

$$\begin{matrix}
M_i\x S^1  & \stackrel{\Phi_k|_{M_i\x S^1}}\to & K(F/F_k,1) \\
                                          
\downarrow_{\text{projection}}  &\hskip-10mm\circlearrowright 
& & & \\
                                          &\nearrow & \\
S^1 &  & \\
\end{matrix}$$



Therefore $\Phi_k$ extends canonically to a continuous map  $\bar\Phi_k:B^5\to K_k$ 
such that $\bar\Phi_k|_{\partial B^5}=\bar\phi_k$. 

This completes the proof of Theorem \ref{Virginia}. 
\qed


\bigbreak

\noindent{\bf Note.} 
 The author could prove that 
if we replace `spun'  with `ribbon' in Lemma \ref{Oklahoma}, 
Lemma \ref{Oklahoma}  
holds.  
See Note \ref{Hartford}.(2)  
and \cite[the paragraph right before Corollary 4.14]{OgasaZ} 
for the definition of ribbon-$(S^2, T^2)$-links.



\begin{note}\label{Hartford}
In \cite{OgasaZ} the author proved the following facts (1)-(3). 
\smallbreak\noindent
(1)   Let $A(t)$ represent 
the 1-st $\Z[t, t^{-1}]$-Alexander polynomial of a $(S^2, T^2)$-link $L$. 
Then 
$
\left.\displaystyle\frac{A(t)}{(t-1)}\right|_{t=1}
=0$
if and only if the alinking number of  $L$ is zero. 

\smallbreak\noindent(2) 
We say that 
a $(S^2, T^2)$-link $L=(K_1,K_2)$ is {\it ribbon} 
if there is an immersion \newline  
$f:B\amalg H\looparrowright S^4$ 
with the following properties, 
where $B$ is a 3-ball and $H$ is an oriented genus one handlebody: 
The self-intersection of $f$ consists of double points and 
is a disjoint union of 2-discs.  
Note that $f^{-1}(\text{each disc})$ is a disjoint union of two 2-discs. 
One of the two disc is included in the interior of $B\amalg H$. 
`(The other disc)$\cap  \partial(B\amalg H)$' is $\partial($the other disc$)$.

Let $A(t)$ represent the 1-st $\Z[t, t^{-1}]$-Alexander polynomial of 
a ribbon-$(S^2, T^2)$-link $L$.     
Then 
$
\begin{vmatrix}
\left.\displaystyle\frac{A(t)}{(t-1)}\right|_{t=1}
\end{vmatrix} 
$
 is the alinking number of  $L$.

\smallbreak\noindent(3)  
Let $A(t)$ and $L$ be as in (2).  
Let $B(t)$ represent the 1-st $\Z[t, t^{-1}]$-Alexander polynomial of 
a ribbon-$(S^2, T^2)$-link $M$.      
If the ribbon-$(S^2, T^2)$-links $L$ and $M$ are cobordant,  
then 
$
\begin{vmatrix}
\left.\displaystyle\frac{A(t)}{(t-1)}\right|_{t=1}
\end{vmatrix}=
\begin{vmatrix}
\left.\displaystyle\frac{B(t)}{(t-1)}\right|_{t=1} 
\end{vmatrix}
$.

The author does not know whether we can remove the condition `ribbon' from (2) \newline and (3). 
\end{note}


\bigbreak
\section{Our new way to investigate surface-link-cobordism}\label{daidai}
\noindent
In order to introduce  our new way, we begin by defining   
terminologies. 

\begin{defn}\label{North Carolina}
Let $n\in\N_{\geqq2}$. 
Let $\mathcal K=(K_1,...,K_n)$ and $\mathcal K'=(K'_1,...,K'_n)$ 
be a surface-link of the standard 4-sphere. 
Let $\{i,j\}\subset\{1,...,n\}.$  
Suppose that $K_i$ is diffeomorphic to $K'_i$. 
If alk$(K_i,K_j)$=alk$(K'_i,K'_j)$ for all distinct $i, j$, 
we say that $\mathcal K$ and $\mathcal K'$ are {\it alinking-equivalent}. 
If there is a surface-link $\mathcal K''=(K''_1,...,K''_n)$ obtained  form $\mathcal K'=(K'_1,...,K'_n)$ 
by changing orders of components such that 
 $\mathcal K$ and $\mathcal K''$ are alinking-equivalent 
and 
such that  for each $i\in\{1,...,n\}$, $K_i$ is diffeomorphic to $K''_i$, 
then we say that 
$\mathcal K$ and $\mathcal K'$ are {\it weakly alinking-equivalent}. 
\end{defn} 

We have the following.

\begin{thm}\label{Ohio}
{\bf (See
\cite[Lemma 4.2]{CassonGordon}
and 
\cite[the first several lines of page 523]{CochranOrr}.)} \newline
$(1)$ 
Let $J$ be a spherical 2-knot $\subset S^4$. 
Let $p$ be a prime power. 
Take the $p$-fold branched cyclic covering space of $S^4$ along $J$, 
and call it $S'$.  
Then $S'$ is a $Z_p$-homology sphere. 

\smallbreak\noindent$(2)$  
Take $S^4\x[0,1]$. 
Let $i\in\{0,1\}$.
Let $J_i$ be a spherical 2-knot contained in $S^4\x\{i\}$. 
Let $p$ be a prime power. 
Take a submanifold $X\subset S^4\x[0,1]$ 
which gives cobordism between $J_0$ and $J_1$. 
$($Note: \cite{Kervaire} ensures that the existence of $X$ $).$  
Take the $p$-fold branched cyclic covering space of $S^4\x[0,1]$ along $X$, 
and call it $M$.  
Let $*\in\Z$.  
Then \newline 
$H_*(M;\Z_p)
\cong 
\begin{cases}
\Z_p&\text{if $*=0,4$}\\
0&\text{else}\\
\end{cases} 
$. 
\end{thm}



\begin{defn}\label{New York}
Let $E=(J,K_1,...,K_m)$ be a surface-link 
of the standard 4-sphere $S^4$.     
Let $n\in\N$. 
Let $M$ be the $n$-fold branched cyclic covering space of $S^4$ along $J$. 
Take the lift  $\widetilde{K}$ of 
the sublink $K=(K_1,...,K_m)$ associated with this branched cyclic covering. 
This submanifold $\widetilde{K}\subset M$ is called 
the {\it$n$-covering-link of $E$ along $J$.}  
\end{defn}

\begin{pr}\label{New Mexico} 
Take $J$, $E$,  $\widetilde{K}$  in Definition \ref{New York}.  
Suppose that $J$ is the trivial spherical 2-knot and that $E$ is semi-boundary.   
Then $\widetilde{K}$ 
is contained in the standard 4-sphere, and 
is an $m\cdot n$-component surface-link.     
\end{pr}

\noindent{\bf Proof of Proposition \ref{New Mexico}.}
Since $J$ is the trivial spherical 2-knot,   
the $n$-fold branched cyclic covering space of $S^4$ along $J$ is the standard 4-sphere.   
Let $N(J)$ be the tubular neighborhood of $J$ in $S^4$. 
Take a $\Z_n$-covering space of $S^4-\text{Int}N(J)$ 
associated with this branched cyclic  covering. 
This $\Z_n$-covering space of $S^4-\text{Int}N(J)$ is regarded as 
the total space of a $\Z_n$-fiber bundle over  $S^4-\text{Int}N(J)$. 
Since any 1-cycle in $K$ is null-homologous in $S^4-\text{Int}N(J)$, 
the restriction of this $\Z_n$-fiber space to $K$ is the trivial bundle. 
Hence Proposition \ref{New Mexico} holds. 
\qed 

\bigbreak
It is convenient to define the following terminologies. 

\begin{defn}\label{Hawaii}  
Let  $*\in\{0,1\}.$
Let $L_*$ be a surface-link contained in the standard 4-sphere $S^4\x\{*\}$.
Let $\mathcal M$ be a connected compact oriented 5-manifold  
which is not necessarily diffeomorphic to $S^4\x[0,1]$. 
Suppose that $\partial \mathcal M$ is 
the disjoint union of two copies of the standard 4-sphere. 
Call one $S^4\x\{0\}$, and the other $S^4\x\{1\}$.  
In Definition \ref{Idaho}, 
replace $S^4\x[0,1]$ 
with 
$\mathcal M$.  
We say that $L_0$ and $L_1$ are {\it cobordant in $\mathcal M$}. 
We say that an embedding map $f$ 
(resp. a submanifold 
$f((F_1\x[0,1])\amalg...\amalg(F_\mu\x[0,1]))\subset\mathcal M$) 
{\it gives cobordism} between $L_0$ and $L_1$ {\it in $\mathcal M$}.  
\end{defn}

\begin{thm}\label{Alaska}   
Let $\iota\in\{0,1\}.$
Let $E_\iota=(K_{\iota0},K_{\iota1},...,K_{\iota,m})$ be a semi-boundary surface-link 
such that $K_{\iota0}$ is the spherical trivial 2-knot.    
Let $*\in\{0,1,...,m\}$. 
Suppose that $E_1$ is cobordant to $E_2$, 
and that 
$C=C_0\amalg C_1\amalg...\amalg C_m\subset S^4\x [0,1]$ 
$($resp. $C_*)$   
gives cobordism between 
$E_0$ and $E_1$ $($resp. $K_{0*}$ and $K_{1*}).$  
Let $p$ be a prime power. 
Take the $p$-fold branched cyclic covering space, $\mathcal M$, of $S^4\x[0,1]$ along $C$. 
Thus we obtain the $p$-covering-link $\mathcal E_\iota$ of $E_\iota$ 
along $K_{\iota0}$. 
Then $\mathcal E_\iota$ is contained in the standard 4-sphere and 
the alinking number associated with the components of $\mathcal E_\iota$ makes sense. 
Then $\mathcal E_0$ is weakly alinking-equivalent to $\mathcal E_1$.  
\end{thm}

\noindent
{\bf Proof of Theorem \ref{Alaska}.} 
Let $\iota\in\{0,1\}.$
Since 
$E_\iota$ is a semi-boundary surface-link and 
$K_{\iota0}$ is the trivial spherical knot, 
$\mathcal E_\iota$ is a $p\cdot m$-component surface-link 
by Proposition \ref{New Mexico}, and let 
$\mathcal E_\iota=(\mathcal E_{\iota1},...,\mathcal E_{\iota,p\cdot m})$.     
Let $\xi\in\{1,...,m\}$. 
Furthermore the lift of $C_\xi$ is a disjoint union of compact connected components 
$C_{\xi1}\amalg,,,\amalg C_{\xi,p}$. 
We can give an order to all elements of $\{C_{\xi,\#}|\xi=1,...,m, \#=1,...,p\}$, 
and call them $A_1,...,A_{p\cdot m}$. 
 Let $\natural\in\{1,...,p\cdot m\}$. 
Let $\partial A_\natural=\mathcal E_{0\natural}\amalg \mathcal E_{1\natural}$. 
Note that  
$\mathcal E_{0\natural}$ is diffeomorphic to $\mathcal E_{1\natural}$ 
and that  
$A_\natural$ is diffeomorphic to $\mathcal E_{0\natural}\x[0,1]$. 
Thus 
$A_1\amalg...\amalg A_{p\cdot m}$ 
(resp $A_\natural$)  
gives cobordism between 
$\mathcal E_0$ and $\mathcal E_1$ 
(resp. $\mathcal E_{0\natural}$ and $\mathcal E_{1\natural}$) in $\mathcal M$.   
Let $\eta\in\Z$.  
By Theorem \ref{New Mexico},  $\mathcal M$ satisfies the condition 
$H_\eta(\mathcal M; \Z_p)
\cong 
\begin{cases}
\Z_p&\text{if $\eta=0,4$}\\
0&\text{else}\\
\end{cases} 
$. 
Since $K_{01}$ and $K_{11}$ are the trivial spherical 2-knot, 
$\partial\mathcal M$ is the disjoint union of two copies of the standard 4-sphere. 
Call one which include $\mathcal E_0$, $S_0$, 
and the other which include $\mathcal E_1$, $S_1$.

\begin{figure}
\includegraphics[width=16cm]{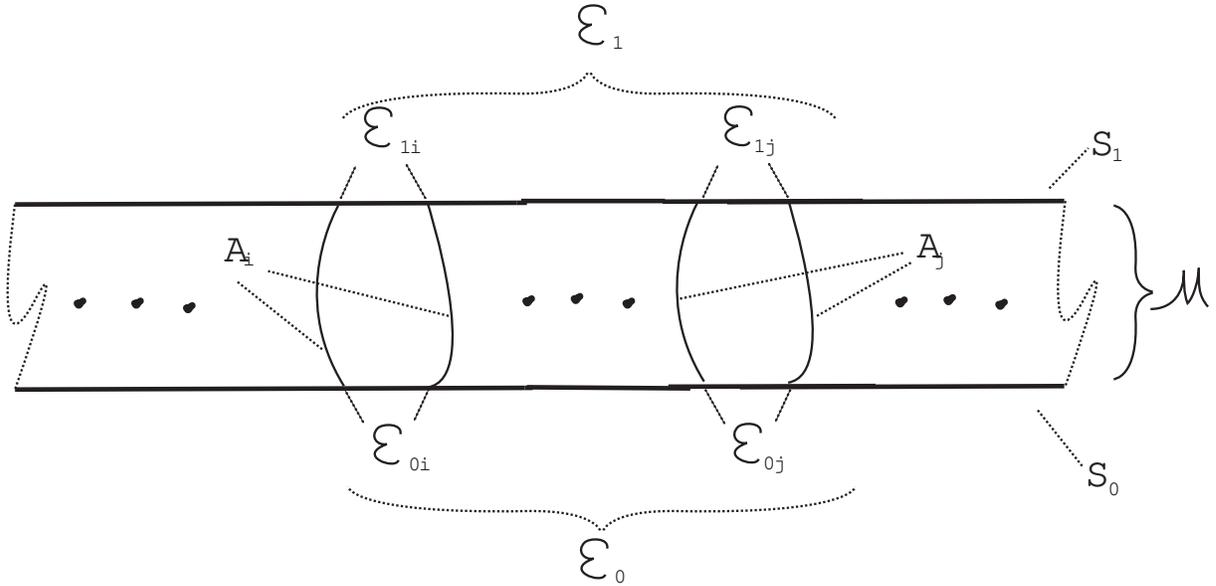}
\vskip-2mm
\caption{\bf  $\mathcal M.$ 
\label{California}}
\vskip-1mm
\end{figure}

Let  $\{i, j\}\subset\{1,...,p\cdot m\}$. 
It suffices to prove that for any distinct $i, j$,   \newline 
the alinking number 
alk$(\mathcal E_{0j}\subset\mathcal E_0, \mathcal E_{0i}\subset\mathcal E_0)$ 
of 
$\mathcal E_{0j}$ in $\mathcal E_0$ around $\mathcal E_{0i}$ in $\mathcal E_0$ 
is equivalent to \newline 
the alinking number 
alk$(\mathcal E_{1j}\subset\mathcal E_1, \mathcal E_{1i}\subset\mathcal E_1)$ 
of 
$\mathcal E_{1j}$ in $\mathcal E_1$ around $\mathcal E_{1i}$ in $\mathcal E_1$.  
Let $\iota\in\{0,1\}$. 
Take a Seifert hypersurface $V_{\iota,i}$ for $\mathcal E_{\iota,i}$.    
Of course $V_{\iota,i}\cap \mathcal E_{\iota,j}\neq\phi$ may hold.  
Note that $V_{0i}\cup A_i\cup V_{1i}$ is a closed oriented 3-manifold 
and a $\Z$-3-cycle in $\mathcal M$. 
Let $\tau\in\Z$. 
Since $H_\tau(\mathcal M;\Z_p)\cong
\begin{cases}
\Z_p&\text{if $\tau=0,4$}\\
0&\text{else}\\
\end{cases}, 
$  
there is a natural number $y_i$ and a $\Z$-4-chain $\Gamma$ 
such that $\partial\Gamma=y_i\cdot(V_{0i}\cup A_i\cup V_{1i})$, 
where $\partial$ denotes the boundary of a chain (not a manifold). 
Furthermore 
we can take an immersion map $g_i$  
from a compact oriented 4-manifold $X_i$ 
to $\mathcal M$ with the following properties: 
$g_i|_{\text{Int} X_i}$ is an embedding map. 
$g_i(\partial X_i)=V_{0i}\cup A_i\cup V_{1i}$ 
$g_i|{\partial X_i}:\partial X_i\to V_{0i}\cup A_i\cup V_{1i}$ 
defines a $Z_{y_i}$-fiber bundle.   ({\it Reason.} Use Thom-Pontrjagin construction.)

Take any embedded circle $C$ in $\mathcal E_{0j}$
Take an embedded circle $C'$ in $\mathcal E_{1j}$ with the following properties: 
There is $Q$ embedded in $A_j$. 
$Q$ is diffeomorphic to $S^1\x[0,1]$. 
$\partial Q=C\amalg C'$, 
where $\partial$ denotes the boundary of a manifold. 
 
Consider $Q\cap X_i$. 
Then $Q\cap X_i$ is a 1-chain. The boundary of the 1-chain is 
algebraically zero points. 
Since 
$A_i\cap A_j=\phi$,  
the points,  which make the boundary of the 1-chain,  are included in 
$V_{0i}\amalg V_{1i}$.   
Hence 
$V_{0i}\cap\mathcal E_{0j}$ and 
$V_{1i}\cap\mathcal E_{1j}$ are algebraically the same number of oriented 0-cells. 
Note that the absolute value of the algebraic number $V_{0i}\cap\mathcal E_{0j}$ 
(resp. $V_{1i}\cap\mathcal E_{1j}$) 
is  
$y_i\cdot$alk$(\mathcal E_{0j}\subset\mathcal E_0, \mathcal E_{0i}\subset\mathcal E_0)$ 
(resp. 
$y_i\cdot$alk$(\mathcal E_{1j}\subset\mathcal E_1, \mathcal E_{1i}\subset\mathcal E_1)$).  
Therefore  \newline
alk$(\mathcal E_{0j}\subset\mathcal E_0, \mathcal E_{0i}\subset\mathcal E_0)=$ 
alk$(\mathcal E_{1j}\subset\mathcal E_1, \mathcal E_{1i}\subset\mathcal E_1)$.  
%
\qed 

\bigbreak
\noindent
{\bf Note.}   
In Theorem \ref{Alaska}, 
if we check the covering-links more elaborately, 
we can prove a stronger condition on
which pair of components of $\mathcal E_0$ and 
which pair of components of $\mathcal E_1$ are equivalent.  
 (The way is written implicitly in this section).  
However we do not need it in order to prove our results in this paper. So we do not discuss it  further here.


\bigbreak
\section{Proof of Main Theorem \ref{Alabama}}\label{Sacramento}
\begin{figure}
\vskip1cm
\includegraphics[width=12cm]{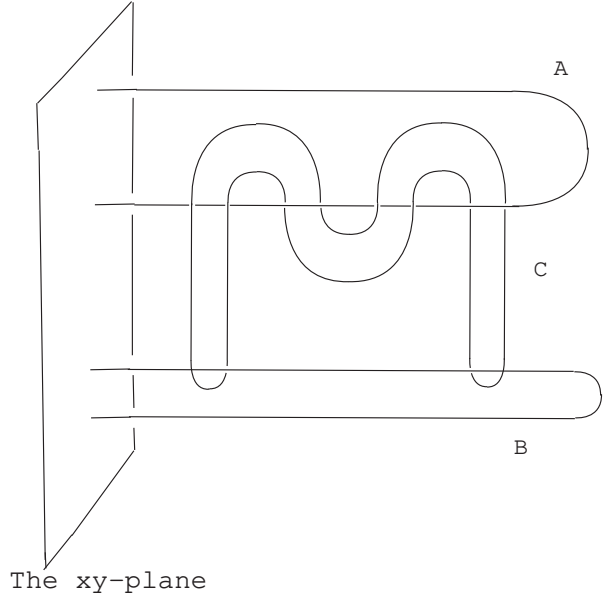}
\caption{\bf  Make a Spun-surface-link from this. 
 \label{Phoenix}}
\end{figure}

\noindent
Make a spun surface-link from $A\amalg B\amalg C$ in Figure \ref{Phoenix}. 
In Figure \ref{Phoenix} we draw the $xy$-plane and omit the $x$-, $y$-, $z$-axes. 
From now we do like this when we draw this kind of figures.  
Let 
$\mathcal L=(P,Q,R)$ be the resultant surface-link such that 
$P$ (resp. $Q$, $R$) is made from $A$ (resp. $B$, $C$). 
Let 
$\mathcal L'=(P',Q',R')$ be the resultant surface-link such that 
$P'$ (resp. $Q'$, $R'$) is made from $B$ (resp. $A$, $C$). 
Thus we have the following. 
\begin{fact}\label{Dover}
There is an orientation-preserving diffeomorphism map 
$f:\R^4\to\R^4$ $($or $S^4\to S^4)$ 
such that 
$f|_{\mathcal L}$ is an order-nonpreserving orientation-preserving diffeomorphism map 
$\mathcal L\to\mathcal L'.$
\end{fact}

\noindent
By the construction, we have the following. 

\begin{fact}\label{Tallahassee}
All sublinks of $\mathcal L$ $($resp. $\mathcal L')$  are the standard links.
$($Recall that we defined that  
$\mathcal L$ $($resp. $\mathcal L')$ itself and the empty set are not sublinks.$)$ 
\end{fact}
 
\noindent
Fact \ref{Tallahassee} implies 
Main Theorem \ref{Alabama}.(1), (2)(ii), (3), (4).(ii), and (5).(ii).

\noindent
Fact \ref{Dover} implies 
Main Theorem \ref{Alabama}.(2)(i),  and (4)(i), and (5).(i). 

\noindent
Fact \ref{Dover} and Theorem \ref{Virginia} imply Main Theorem \ref{Alabama}.(5).(ii). 

We prove Main Theorem \ref{Alabama}.(6). 
Let $p$ be a sufficiently large prime power. 
Take the $p$-covering-link $\mathcal Z$ (resp. $\mathcal Z'$) 
of $\mathcal L$ (resp. $\mathcal L'$) along $P$ (resp. $P'$). 

By the definition of  $\mathcal L$, we have the following: 
Let $\mathcal Z_i$ and  $\mathcal Z_j$ be different components of  $\mathcal Z$.  
The alinking number of  $\mathcal Z_i$ around  $\mathcal Z_j$ is one or zero. 

By the definition of  $\mathcal L'$, we have the following: 
There are two different components $\mathcal Z'_k$ and  $\mathcal Z'_l$ of  $\mathcal Z'$
such that the alinking number of  $\mathcal Z'_k$ around  $\mathcal Z'_l$ is two. 

Therefore $\mathcal Z$ and $\mathcal Z'$ are not weakly alinking-equivalent. 
By Theorem \ref{Alaska}, 
$\mathcal L$ and $\mathcal L'$ are not cobordant. 

This completes the proof of Main Theorem \ref{Alabama}. 
\qed

\bigbreak
\section{A theorem}\label{New Jersey}

\noindent
We give another partial solution to Problem \ref{koremo}. 
Theorem \ref{standard} lets us naturally formulate the following problem, 
which is included in Problem \ref{koremo}. 

\begin{prob}\label{sui}   
 Find a new $(S^2, T^2)$-link which is not cobordant to the standard link. 
\end{prob}



Theorem \ref{main} is a partial solution to this problem 
since $\mathcal L^{(0)}$ in Proof of  Theorem \ref{main} is the standard link. 

By Theorem \ref{standard} we have the following: 
if we replace `the standard link' with `a boundary link' in 
Problem \ref{sui},   
the new one is the same as the old one.    
By Corollary \ref{coro} we have the following:  
if we replace `$(S^2,T^2)$-link' with `semi-boundary $(S^2,T^2)$-link' in 
Problem \ref{sui},   
the new one is the same as the old one.    
%
%
%
If we replace $(S^2, T^2)$ with $(S^2, S^2)$ in Problem \ref{sui},  
it is Problem \ref{outstanding}
\footnote{
Note the following facts 
written in \cite{CochranOrr}: 
A natural `even dimensional spherical link' version 
of Problem \ref{outstanding} is 
`Are all  even-dimensional spherical links slice?'
If $n\geqq3$, a natural `spherical $n$-link' version 
of Problem \ref{outstanding} is 
`Are all  spherical $n$-links cobordant to boundary links?'
Note that the latter problem includes the former one.   
The latter one and Theorem \ref{standard} let us naturally regard Problem \ref{sui} as a surface-link version of Problem \ref{outstanding}.  
%
}. 
%
%


\begin{thm}\label{main} 
There is a set 
$\{(S^2, T^2)$-links $\mathcal L^{(i)}=(\mathcal J^{(i)}, \mathcal K^{(i)})| i,j\in\N\cup\{0\}$,  
$\mathcal L^{(i)}$ is  non-cobordant to $\mathcal L^{(j)}$ 
for any distinct $i,j\}$  
with the following properties: 

\begin{enumerate}
%
\item
The alinking number of $\mathcal L^{(i)}$ is the same as that of $\mathcal L^{(j)}$ 
for any distinct $i,j$. 

\item
The Sato-Levine invariant of  $\mathcal L^{(i)}$ is the same as that of $\mathcal L^{(j)}$ 
for any distinct $i,j$.

\item
The Cochran sequence of  $\mathcal L^{(i)}$ is the same as that of $\mathcal L^{(j)}$ 
for any distinct $i,j$. 

\item
If there is an element  $*\in\N_{\geqq2}\cup\{\omega\}$ such that
for any distinct $i,j$,   
 we can define the Orr invariant  
$\theta_*(\mathcal L^{(i)}, \tau^{(i)})$ 
for a meridian $\tau^{(i)}$ 
and 
 $\theta_*(\mathcal L^{(j)}, \tau^{(j)})$  
for a meridian $\tau^{(j)}$,      
then $\theta_*(\mathcal L^{(i)})=\theta_*(\mathcal L^{(j)})=0$. 

\item
Our new way in \S\ref{daidai} implies that   
 $\mathcal L^{(i)}$ is non-cobordant to $\mathcal L^{(j)}$ 
for any distinct $i,j$.  
%
%
%
\end{enumerate}
\end{thm}



%




\begin{figure}
\includegraphics[width=14cm]{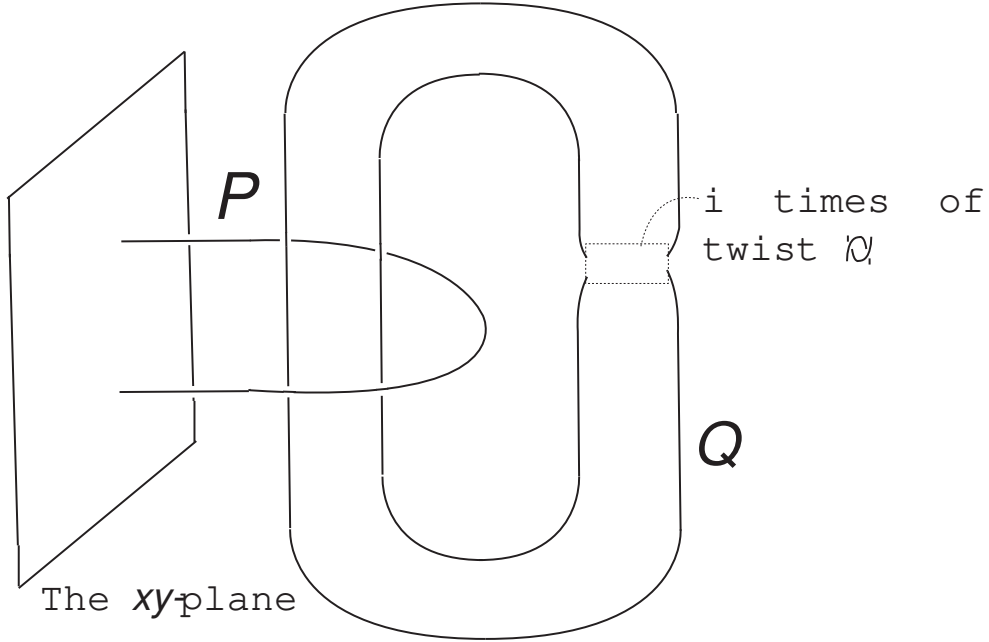}
\vskip-14mm
\caption{\bf  
$P\amalg Q$ (resp. $P$, $Q$) 
is made into 
the spun $(S^2, T^2)$-link 
$\mathcal L^{(i)}=(\mathcal J^{(i)}, \mathcal K^{(i)})$  
(resp. $\mathcal J^{(i)}$, $\mathcal K^{(i)}$) by the rotation.   
 \label{Little Rock}}
\end{figure}

\begin{figure}
\includegraphics[width=10cm]{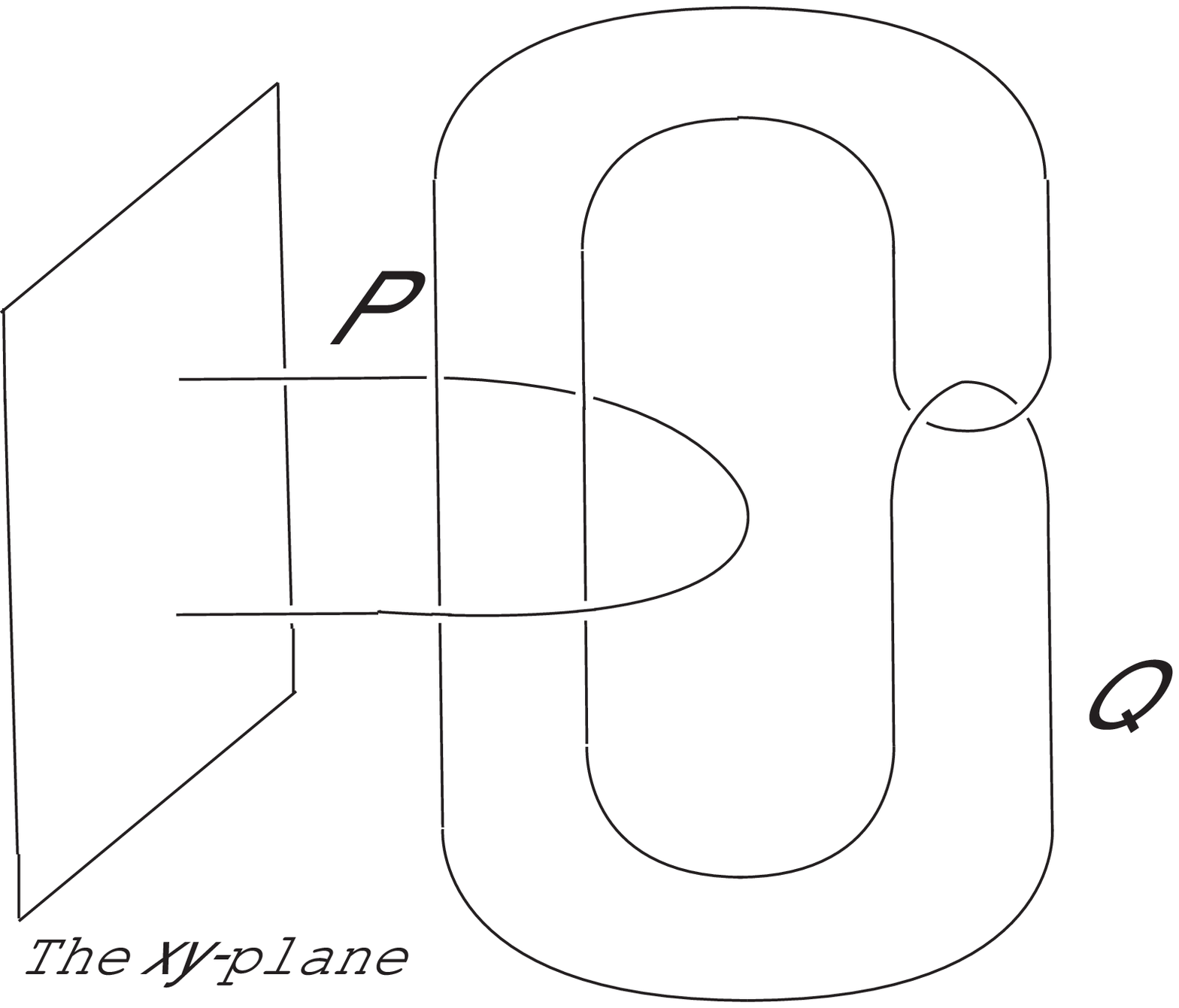}

\vskip-17mm
%
\caption{\bf  
The $i=1$ case of Figure \ref{Little Rock}, that is, $\mathcal L^{(1)}$
\label{Colorado}}
\end{figure}

\noindent{\bf Proof of Theorem \ref{main}.}
Make a spun $(S^2, T^2)$-link 
$\mathcal L^{(i)}=(\mathcal J^{(i)}, \mathcal K^{(i)})$  
(resp. $\mathcal J^{(i)}$, $\mathcal K^{(i)}$) 
from 
$P\amalg Q$ (resp. $P$, $Q$) in Figure \ref{Little Rock} 
by the rotation around the $xy$-plane as the axis.   
Figure \ref{Colorado} is the $i=1$ case of Figure \ref{Little Rock}. 
We prove that 
$\{\mathcal L^{(i)}| i\in\N\cup\{0\}\}$ satisfies the condition (1)-(5) of Theorem \ref{main}. 
By the construction $\mathcal J^{(i)}$ is the trivial spherical 2-knot for each $i$. 

\begin{figure}
\vskip5mm
\includegraphics[width=12cm]{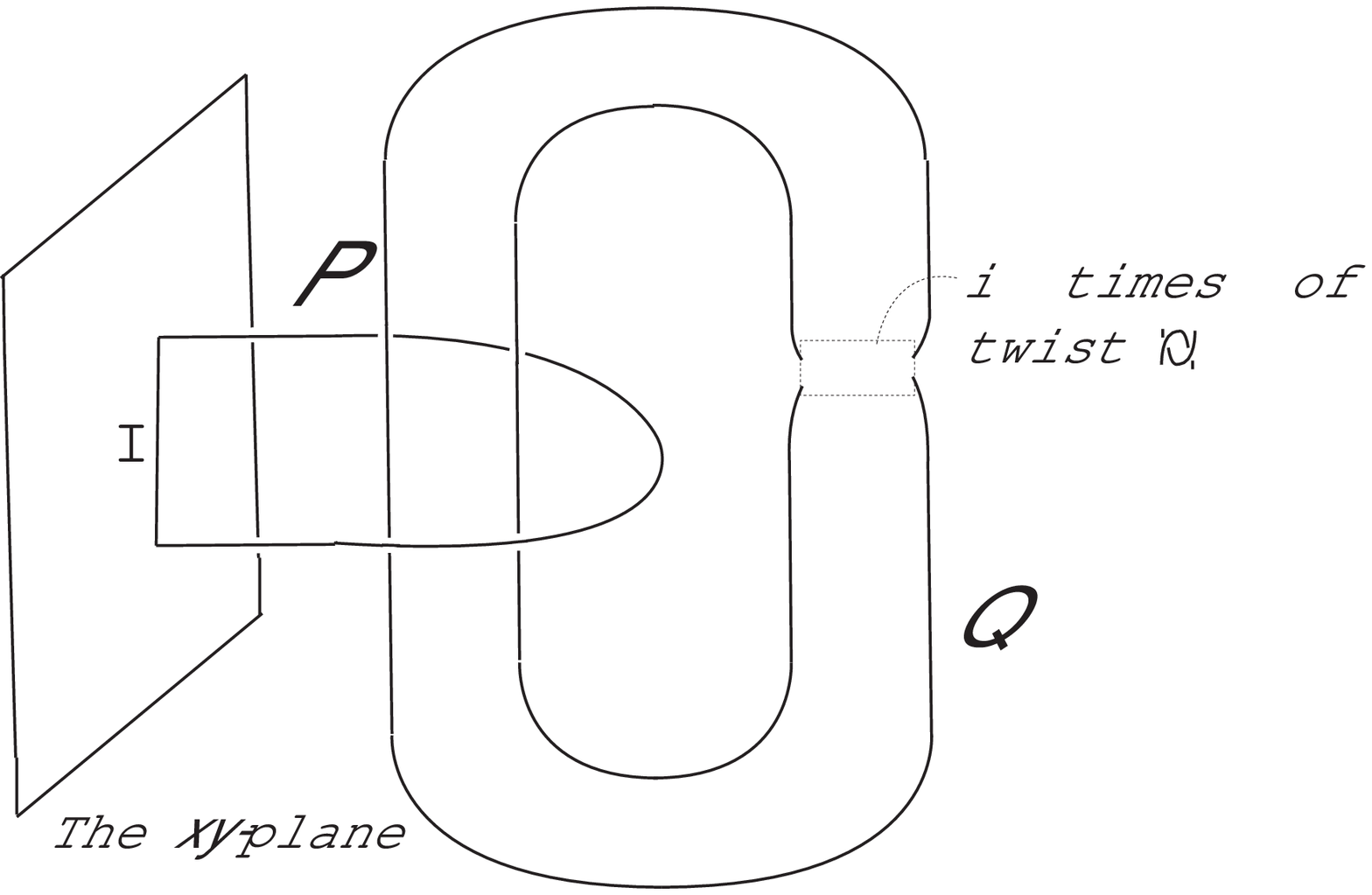}   
\vskip-14mm
\caption{\bf  
We add $I$ to Figure \ref{Little Rock}   
 \label{Connecticut}}
\end{figure}

\bigbreak
We prove that $\{\mathcal L^{(i)}| i\in\N\cup\{0\}\}$ satisfies Theorem \ref{main}.(1).  
Take a 1-link $(P\cup I, Q)$ in $\{(x,y,z)|z\geqq0\}$ as in 
Figure \ref{Connecticut}.  
Note that there is a Seifert surface $V$ (resp. $W$) for $P\cup I$ (resp. $Q$) 
such that $V\cap Q=W\cap(P\cup I)=\phi$. 
Note that 
$V\cup W\subset\{(x,y,z)|z\geqq0\}$. 
We can suppose that 
$V\cap\{(x,y,z)|z=0\}=I$ 
and 
$W\cap\{(x,y,z)|z=0\}=\phi$. 
When we rotate $P\amalg Q$ and make $\mathcal L^{(i)}=(\mathcal J^{(i)}, \mathcal K^{(i)})$, 
rotate $V$ (resp. $W$) together. 
Then the result  
$\mathcal V^{(i)}$ (resp. $\mathcal W^{(i)}$) is a Seifert hypersurface for 
$\mathcal J^{(i)}$ (resp. $\mathcal K^{(i)}$)   
such that   
$\mathcal V^{(i)}\cap \mathcal K^{(i)}=\mathcal W^{(i)}\cap\mathcal J^{(i)}=\phi$.  
The existence of $\mathcal V^{(i)}$ implies 
that $\{\mathcal L^{(i)}| i\in\N\cup\{0\}\}$ satisfies Theorem \ref{main}.(1).  
We do not use $\mathcal W^{(i)}$ here
but we will use it from now.

\begin{figure}
\bigbreak
\includegraphics[width=90mm]{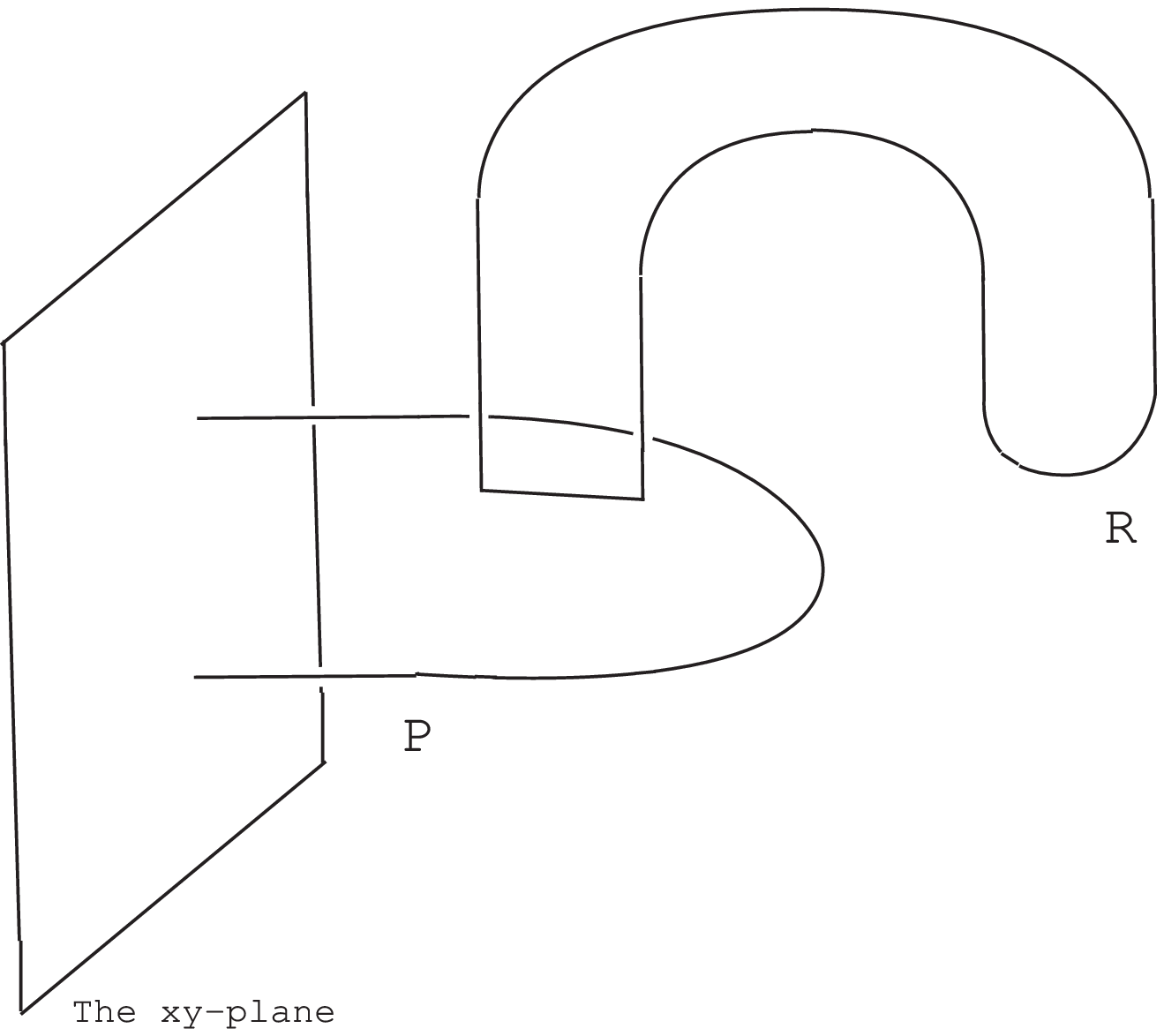}
\caption{\bf  
$P\amalg R$ is a part of a 1-link whose spun link is the derivative of 
$\mathcal L^{(i)}$.   \label{ZDelaware}}
%
%
\end{figure}

\bigbreak
We prove that $\{\mathcal L^{(i)}| i\in\N\cup\{0\}\}$ satisfies Theorem \ref{main}.(2).  
By Theorem \ref{main}.(1),  
we can define the Sato-Levine invariant for $\mathcal L^{(i)}$. 
We can suppose that $V$ intersects $W$ transversely 
and that $V\cap W$ is $R$ in Figure \ref{ZDelaware}.  
We can suppose that 
$\mathcal V^{(i)}$ intersect $\mathcal W^{(i)}$ transversely 
and that 
$\mathcal V^{(i)}\cap\mathcal W^{(i)}$ is 
the result $\mathcal T^{(i)}$ of rotating $R$ in Figure  \ref{ZDelaware}  
around the $xy$ plane 
when we rotate $P\amalg Q$ and make 
$\mathcal L^{(i)}=(\mathcal J^{(i)}, \mathcal K^{(i)})$. 
Note that $\mathcal T^{(i)}$ is diffeomorphic to the torus. 
Take a point in $R$. The point becomes a circle when we make $R$ into $\mathcal T^{(i)}$ by the rotation. 
Let $\sigma$ be the spin structure on the circle 
which is defined by $\mathcal V^{(i)}$ and  $\mathcal W^{(i)}$.  
Then  $[S^1, \sigma]=0\in\Omega^{\rm Spin}_1$.   
Hence we have the following: 
Let $\tau$ be the spin structure on  $\mathcal T^{(i)}$ 
which is defined by $\mathcal V^{(i)}$ and  $\mathcal W^{(i)}$.  
Then $[T^2, \tau]=0\in\Omega^{\rm Spin}_2$.  
See \cite[\S5.6]{Gompf} and \cite[\S IV]{Kirby} for the spin structure.  
%
Since the Sato-Levine invariant of $\mathcal L^{(i)}$ is 
$[T^2, \tau]\in\Omega^{\rm Spin}_2\cong\Z_2$ (see \cite{OgasaSL}), 
 the Sato-Levine invariant of $\mathcal L^{(i)}$ is zero. 
Hence $\{\mathcal L^{(i)}| i\in\N\cup\{0\}\}$ satisfies Theorem \ref{main}.(2).

\bigbreak
We prove that $\{\mathcal L^{(i)}| i\in\N\cup\{0\}\}$ satisfies Theorem \ref{main}.(3). 
By Theorem \ref{main}.(1), 
we can define the Cochran sequence for $\mathcal L^{(i)}$.
The derivative $D(\mathcal L^{(i)})$ of $\mathcal L^{(i)}$ is 
the spun $(S^2, T^2)$-link of Figure  \ref{ZDelaware}.  
It is the standard link.  
Hence $\beta^j(\mathcal L^{(i)})=0 $ if $j\geqq2$.   
By this fact and Theorem \ref{main}.(2),  we have $\{\beta^j(\mathcal L^{(i)})\}_{j\in\N}$ is trivial. 
Hence $\{\mathcal L^{(i)}| i\in\N\cup\{0\}\}$ satisfies Theorem \ref{main}.(3).

\bigbreak
By Theorem \ref{Virginia}, 
$\{\mathcal L^{(i)}| i\in\N\cup\{0\}\}$ satisfies Theorem \ref{main}.(4).

\bigbreak
We prove that $\{\mathcal L^{(i)}| i\in\N\cup\{0\}\}$ 
satisfies Theorem \ref{main}.(5).
\begin{figure}
\includegraphics[width=10cm]{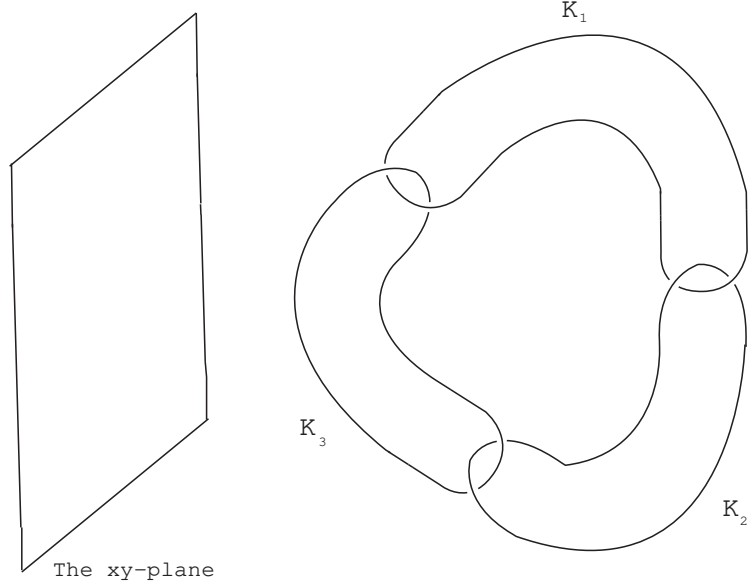}
\caption{\bf  
A 1-link $(K_1,K_2,K_3)$ whose spun $(S^2, T^2)$-link is 
the 3-covering-link $\mathcal L^{(1,3)}
=(\mathcal K_1^{(1,3)},\mathcal K_2^{(1,3)},\mathcal K_3^{(1,3)})$  
of $\mathcal L^{(1)}$.  $K_i$ is made into $\mathcal K_i^{(1,3)}$ by the rotation.  
\label{ZFlorida}}
\end{figure}
%
%
%
%
%
%
%
Let $p$ be a prime power. 
By Proposition \ref{New Mexico}, 
the $p$-covering-link $\mathcal L^{(i,p)}$ of $\mathcal L^{(p)}$ is  
a $\underbrace{(T^2,...,T^2)}_{p}$-link $(\mathcal K^{(i,p)}_1,...,\mathcal K^{(i,p)}_p)$.  
Figure \ref{ZFlorida} 
is the case $i=1$ and $p=3$. 
Let $\nu\in\{1,2,3\}$. 
Note that $K_\nu$ in Figure \ref{ZFlorida} 
is made into  $\mathcal K^{(i,p)}_\nu$. 
%
%

By the construction of $\mathcal L^{(i,p)}$, we have the following: 
If $p$ is sufficient large,  
then for any distinct $\alpha, \beta\in\{1,...,p\}$, 
the alinking number of 
  $\mathcal K^{(i,p)}_\alpha$ around $\mathcal K^{(i,p)}_\beta$ 
is $i$ or 0. 
By Theorem \ref{Alaska},  
$\mathcal L^{(i)}$ is not cobordant to $\mathcal L^{(j)}$ for any distinct $i,j\in\N\cup\{0\}$. 
Hence 

\noindent 
$\{\mathcal L^{(i)}| i\in\N\cup\{0\}\}$ satisfies Theorem \ref{main}.(5).

This completes the proof of Theorem \ref{main}.
\qed

\bigbreak\noindent{\bf Note.} 
(1) By using Theorem \ref{Virginia}, 
the author could make infinitely many examples like 
the set $\{\mathcal L^{i}\}$ in Theorem \ref{main}. 

\smallbreak
\noindent
(2)The author could prove that 
we can define the Cochran sequence for 
the spun $(S^2, T^2)$-link of any semi-boundary 2-component 1-link, 
where \newline
(any component of the 1-link)$\cap$(the axis $\R^3$) is $\phi$ or the interval, 
and that it is trivial.  In order to prove this fact, 
we use the $T^2_{bd}$ spin structure (see \cite{Kirby})  
as in the proof of Theorem \ref{main}.   

\bigbreak
\section{Discussion}\label{Florida}
\noindent 
Let $i\in\N\cup\{0\}$. 
Take the $(S^2, T^2)$-link $\mathcal L^{(i)}=(\mathcal J^{(i)}, \mathcal K^{(i)})
\subset S^4$ in Proof of Theorem \ref{main}.    
Suppose that $\mathcal L^{(i)}\subset B^4\subset S^4$. 
Take the standard $T^2$-knot $G\subset S^4-B^4$. 
By using an embedded 3-dimensional 1-handle $\subset S^4$, 
connect $\mathcal J^{(i)}$ and $G$. 
Thus we obtain a  
$(T^2, T^2)$-link $\mathcal M^{(i)}=(\mathcal E^{(i)}, \mathcal F^{(i)})$ 
by this connected-sum. 
Let $\mathcal {M'}^{(i)}=(\mathcal {E'}^{(i)}, \mathcal {F'}^{(i)})$  be 
a $(T^2, T^2)$-link which is made from $\mathcal M^{(i)}$ by changing the order. 
We ask a question. 

\begin{que}\label{Denver}
Let $i\in\N$. 
Are the above 
 $\mathcal M^{(i)}$  
and 
$\mathcal {M'}^{(i)}$ 
cobordant?
\end{que}

Let $i\in\N$. 
By Theorem \ref{main} and the construction of 
 $\mathcal M^{(i)}$ 
and 
$\mathcal {M'}^{(i)}$, 
we have the following: 

\smallbreak\noindent$(1)$
%
%
alk$(\mathcal {E}^{(i)}\subset \mathcal {M}^{(i)}, \mathcal {F}^{(i)}\subset \mathcal {M}^{(i)})$
=alk$(\mathcal {E'}^{(i)}\subset \mathcal {M'}^{(i)}, \mathcal {F'}^{(i)}\subset \mathcal {M'}^{(i)})$

\hskip2mm
alk$(\mathcal {F}^{(i)}\subset \mathcal {M}^{(i)}, \mathcal {E}^{(i)}\subset \mathcal {M}^{(i)})$
=alk$(\mathcal {F'}^{(i)}\subset \mathcal {M'}^{(i)}, \mathcal {E'}^{(i)}\subset \mathcal {M'}^{(i)})$

\smallbreak\noindent$(2)$
The 2-component-Sato-Levine invariant of $\mathcal M^{(i)}$ 
is the same as that of $\mathcal {M'}^{(i)}$.

\smallbreak\noindent$(3)$  
The Cochran sequence of $\mathcal M^{(i)}$ 
is the same as that of $\mathcal {M'}^{(i)}$. 

\smallbreak\noindent$(4)$ 
$\pi_{\mathcal M^{(i)}}$ and $\pi_{\mathcal {M'}^{(i)}}$ are Stallings-equivalent.

\smallbreak\noindent$(5)$
Let  
$\#\in\N_{\geqq2}\cup\{\omega\}$. 
The following two are equivalent. 

(I) 
We can define the Orr invariant 
$\theta_\#(\mathcal {M}^{(i)}, \tau)$ for a meridian $\tau$

(II)
We can define  
$\theta_\#(\mathcal {M'}^{(i)}, \tau')$ for a meridian $\tau'$.

Suppose that (I) (resp. (II)) holds. 
Then  
$\theta_*(\mathcal {M}^{(i)})=\theta_*(\mathcal {M'}^{(i)})=0$.

\bigbreak
However we cannot use the alinking number 
associated with covering-links in general 
as we can in Theorem \ref{Alaska}. 
{\it Reason.} The covering space of $S^4$ along  
$\mathcal {E'}^{(i)}$ (resp. $\mathcal {F'}^{(i)}, \mathcal {E'}^{(i)}, \mathcal {F'}^{(i)}$) 
is not 
an integral homology sphere.

Even if 
we can define the alinking number associated with covering-links 
in some special cases, 
it may not be a cobordism invariant.  
{\it Reason}:
We omit $(i)$ in 
$\mathcal {M}^{(i)}=(\mathcal {E}^{(i)}, \mathcal {F}^{(i)})$ 
and 
$\mathcal {M'}^{(i)}=(\mathcal {E'}^{(i)}, \mathcal {F'}^{(i)})$ 
from now for the convenience. 
Suppose that there is a submanifold $E\amalg F\subset S^4$ 
which gives cobordism between $\mathcal M$ and  $\mathcal {M'}$. 
Note that  $E$ (resp. $F$) is diffeomorphic to $T^2\x [0,1]$. 
Suppose that 
$E$ (resp. $F$) 
gives cobordism between 
$\mathcal E$ and  $\mathcal {E'}$
(resp. $\mathcal F$ and  $\mathcal {F'}$). 
%
%
%
%
The covering space $X$ of $S^4\x[0,1]$ along $E$ (resp. $F$) does not satisfy 
the condition $H_i(X ,\Z)\cong H_i(S^4,\Z)$ for all $i$.  
The way in Proof of Theorem \ref{Alaska} cannot be used.

\bigbreak
Recall that the above (4) follows from \cite[5.2 Theorem]{Stallings}.  
By the way, \cite[5.2 Theorem]{Stallings} follows from \cite[5.1 Theorem]{Stallings}.  
We may be able to prove that  
 $\mathcal M$ 
and 
$\mathcal M'$  
are non-cobordant 
if we use \cite[5.1 Theorem]{Stallings} as follows. 
Note that 
the tubular neighborhood of 
$\mathcal M$ 
in $S^4$ is diffeomorphic to 
$\mathcal M\x D^2$,  
and call it $\mathcal M\x D^2$. 
By \cite[5.1 Theorem]{Stallings}, for each $k\in\N$, natural inclusion maps induce the following isomorphisms. 
%
%

$$\begin{matrix}
\pi_1(S^4-\text{Int}\mathcal M\x D^2)/(\pi_1(S^4-\text{Int}\mathcal M\x D^2))_k
  \\
\downarrow_{\alpha, \cong}   \\
\pi_1(S^4\x[0,1]-\text{Int}(E\amalg F)\x D^2)/(\pi_1(S^4\x[0,1]-\text{Int}(E\amalg F)\x D^2))_k
                             \\
\end{matrix}$$

\bigbreak
$$\begin{matrix}
\pi_1(S^4-\text{Int}\mathcal M'\x D^2)/(\pi_1(S^4-\text{Int}\mathcal M'\x D^2))_k
  \\
\downarrow_{\beta, \cong}   \\
\pi_1(S^4\x[0,1]-\text{Int}(E\amalg F)\x D^2)/(\pi_1(S^4\x[0,1]-\text{Int}(E\amalg F)\x D^2))_k
                             \\
\end{matrix}$$

\noindent
These give an isomorphism.

$$\begin{matrix}
\pi_1(S^4-\text{Int}\mathcal M\x D^2)/(\pi_1(S^4-\text{Int}\mathcal M\x D^2))_k
  \\
\downarrow_{\gamma,\cong}   \\
\pi_1(S^4-\text{Int}\mathcal M'\x D^2)/(\pi_1(S^4-\text{Int}\mathcal M'\x D^2))_k
  \\
\end{matrix}$$


\noindent
By 
the construction of $\mathcal M$  and $\mathcal {M'}$, 
there is an orientation preserving diffeomorphism 
$f:S^4\to S^4$ 
such that 
$f_{\mathcal M}$ is an orientation-preserving, order-nonpreserving diffeomorphism 
$\mathcal M\to\mathcal M'$. 
This gives another isomorphism.

$$
\pi_1(S^4-\text{Int}\mathcal M\x D^2)
\stackrel{\delta,\cong}\to 
\pi_1(S^4-\text{Int}\mathcal M'\x D^2) 
$$





The combination of these isomorphisms, $\alpha$, $\beta$, $\gamma$, and $\delta$ 
may give a restriction to a link-type of $\mathcal M$.

By using the above way, 
we may prove that 
$\mathcal L$ and $\mathcal L'$ in Main theorem \ref{Alabama}  
are noncobordant without using Main theorem \ref{Alabama}(6). 

\smallbreak
By using \cite[5.2 Theorem]{Stallings},  
for a pair $i,j$,  
we may prove that $\mathcal L^{(i)}$ and $\mathcal {L'}^{(i)}$ in Theorem \ref{main} 
are noncobordant  without using Theorem \ref{main} (5). 
(Note. It may be proved as follows. 
$\mathcal L^{(i)}$ and $\mathcal {L}^{(j)}$ may not be Stallings-equivalent. 
Use this fact.)

However, if we can do these, they must be more complicated than the way in this paper.

\bigbreak
We give another cobordism invariant. 
We use the notations in Theorem \ref{Alaska} and Proof of Theorem \ref{Alaska}. 
Take the cobordant surface-links $E_1$ and $E_2$. 
Let $r$ be a sufficiently large natural number. 
Let $p$ be $2^r$. 
For any distinct $i,j$, take $\mathcal E_{0i}$ and $\mathcal E_{0j}$. 
By Theorem \ref{alkcob},  
the following two are equivalent. 

(I)  
The alinking number 
$\mathcal E_{0i}$ (resp. $\mathcal E_{0j}$) around $\mathcal E_{0j}$ (resp. $\mathcal E_{0i}$)
is zero.   

(II)  
The alinking number 
$\mathcal E_{1i}$ (resp. $\mathcal E_{1j}$) around $\mathcal E_{1j}$ (resp. $\mathcal E_{1i}$)
is zero.   

\noindent 
We suppose (I) (resp. (II)).  
Let $*\in\{0,1\}.$  
Hence we can suppose that 
$V_{*i}$ and $V_{*j}$ are special Seifert hypersurface. 
Let $F_*$ be a closed oriented 2-dimensional submanifold $V_{*i}\cap V_{*j}$.   
We induce the spin structure $\phi_*$ on $F_*$ 
by using $V_{*i}$, and $V_{*j}$, and $S^4_*$. 
The spin cobordism class 
$[(F_*, \phi_*)]\in\Omega^{\text{spin}}_2\cong\Z_2$ 
is the 2-component Sato-Levine invariant of  
the 2-component link 
$(\mathcal E_{*i}, \mathcal E_{*j})$.   

\begin{cla}\label{Honolulu}  
$[(F_0, \phi_0)]=[(F_1, \phi_1)]$.  
\end{cla}

\noindent
{\bf Proof of Claim \ref{Honolulu}.}
Since $p$ is an even number, 
$y_i$ and $y_j$ are odd numbers. 
By the existence of $X_i$ and $X_j$, we have  
$y_i\cdot y_j\cdot[(F_0, \phi_0)]=y_i\cdot y_j\cdot[(F_1, \phi_1)]$) $\in\Omega^{\text{spin}}_2$. 
Since $y_i$ is and $y_j$ are odd numbers, 
$[(F_0, \phi_0)]=[(F_1, \phi_1)]$) $\in\Omega^{\text{spin}}_2$. 
This completes the proof. \qed
\bigbreak

Let $*\in\{0,1\}.$  
Under the above condition,  
 for any distinct $i,j,k$,  
make a 3-component sublink 
$Z=(\mathcal E_{*i}, \mathcal E_{*j}, \mathcal E_{*k})$  
of $\mathcal E_*$.  
Let $\#=i,j,k$. 
Let $V_{*\#}$ be a special Seifert hypersurface associated with $Z$. 
Let $C_*$ be the 1-dimensional manifold $V_{*i}\cap V_{*j}\cap V_{*k}$.  
We induce the spin structure $\tau_*$ on $C_*$ 
by using 
$V_{*i}$, $V_{*j}$,  $V_{*k}$, and $S^4_*$. 
The spin cobordism class 
$[(C_*, \tau_*)]\in\Omega^{\text{spin}}_1\cong\Z_2$ 
is the 3-component Sato-Levine invariant of  
the 3-component link 
$(\mathcal E_{*i}, \mathcal E_{*j}, \mathcal E_{*k})$.

\begin{cla}\label{Boise}  
$[(C_0, \tau_0)]=[(C_1, \tau_1)]$. 
\end{cla}

\noindent
{\bf Proof of Claim \ref{Boise}.}
By the existence of $X_i$, $X_j$ and $X_k$, \newline
$y_i\cdot y_j\cdot y_k\cdot [(C_0, \tau_0)]$
$=y_i\cdot y_j\cdot y_k\cdot [(C_1, \tau_1)]\in\Omega^{\text{spin}}_1$. 
Since $y_i$, $y_j$ and $y_k$ are odd numbers, 
$[(C_0, \tau_0)]=[(C_1, \tau_1)]$. \qed

\begin{figure}
\includegraphics[width=7cm]{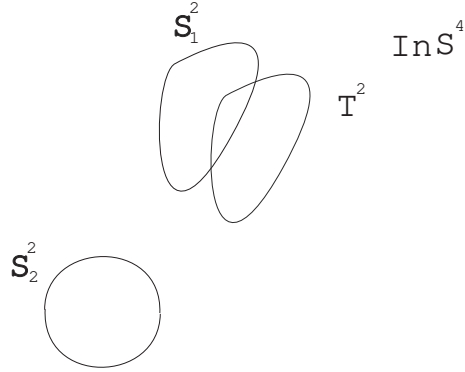}

\caption{\bf  
$(S^2_2,T^2)$ and  $S^2_1$
 \label{Montgomery}}
%
\end{figure}

\begin{figure}
\vskip-15mm
\includegraphics[width=100mm]{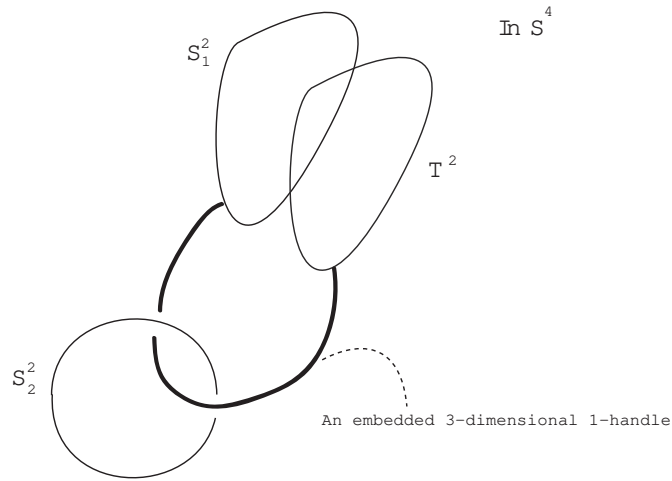}

\caption{\bf  
$(S^2_1\#T^2, S^2_2)$
 \label{Juneau}}
\end{figure}

\bigbreak
We show an example of two noncobordant surface-links 
such that we cannot distinguish the nonequivalence of their cobordism classes   
by using the alinking number associated with a covering-link in the way of Theorem \ref{Alaska}, 
but that we can do by using the 2-component Sato-Levine invariant as above. 
Let the standard $(S^2, T^2)$-link be the one of the two surface links. 
See Figure \ref{Montgomery}.  
Take a $(S^2, T^2)$-link $(S^2_1,T^2) \subset B^4\subset S^4$ 
whose 2-component-Sato-Levine invariant is nonzero. 
Take the trivial spherical 2-knot $S^2_2\subset S^4-B^4$. 
See Figure \ref{Juneau}.   
Make a band-sum $S^2_1\#T^2$ of $S^2_1$ and $T^2$ 
by using an embedded 3-dimensional 1-handle 
as drawn in Figure \ref{Juneau} 
so that one of the above $[(F_*, \phi_*)]$ is nontrivial.  
Then we obtain a $(S^2, T^2)$-link $(S^2_2, S^2_1\#T^2)$. 
Let this link the other of the two.

We show an example of two noncobordant surface-links 
such that we cannot distinguish the nonequivalence of their cobordism classes   
by using the alinking number associated with a covering link in the way of Theorem \ref{Alaska} 
or by using the 2-component Sato-Levine invariant as above, 
but that we can do by using the 3-component Sato-Levine invariant as above. 
Take a 3-component surface-link $\subset B^4\subset S^4$ 
whose 3-component(resp. 2-component)-Sato-Levine invariant is nontrivial (resp. trivial). 
Take the trivial 2-knot $\subset S^4-B^4$. 
Make a band-sum by using two (not one) embedded 3-dimensional 1-handles 
in a similar way to the above one,   
and obtain a 2-component surface-link $L$. 
Let $L$ be the one of the two. 
Let the other of the two  
the standard surface link which is orientation-preserving diffeomorphic to $L$.

\bigbreak
Investigate relations between 
these invariants defined in this section 
and 
the other invariants defined or cited in the previous sections.

\bigbreak

\bigbreak\noindent
Eiji Ogasa:  Computer Science, Meijigakuin University, Yokohama, Kanagawa, 244-8539, Japan 
\quad pqr100pqr100@yahoo.co.jp  \quad
ogasa@mail1.meijigkakuin.ac.jp 

\end{document}